\documentclass[reqno]{amsart}
\usepackage{amssymb}
\usepackage{amsfonts}
\usepackage{amsmath}
\usepackage{amsthm, amsmath, amssymb, enumerate}
\usepackage{amssymb, enumitem}
\usepackage{bbm}
\usepackage{amscd}
\usepackage[all]{xy}
\usepackage[hidelinks]{hyperref}
\usepackage{amsmath}
\usepackage{amssymb}
\usepackage{amsthm}
\usepackage[left=2cm,right=2cm,bottom=3cm,top=3cm]{geometry}
\usepackage[sort&compress,numbers]{natbib}
\usepackage{tikz}
\usepackage{tikz-cd}
\usepackage{wrapfig}
\usepackage{comment}
\usepackage{centernot}
\usepackage{mathtools}
\usepackage{stmaryrd}

\newtheorem{theorem}{Theorem}[section]
\newtheorem{proposition}[theorem]{Proposition}
\theoremstyle{definition}
\newtheorem{lemma}[theorem]{Lemma}
\newtheorem{claim*}[theorem]{claim*}
\newtheorem{definition}[theorem]{Definition}
\newtheorem{corollary}[theorem]{Corollary}

\newtheorem{remark}[theorem]{Remark}

\AtBeginDocument{   \def\MR#1{}}

\begin{document}
\title{A general Universal Coefficient Theorem, and applications}
\author{Luigi Caputi}
\address{Dipartimento di Matematica, Universit\`{a} di Bologna, Piazza di
Porta S. Donato, 5, 40126 Bologna,\ Italy}
\email{luigi.caputi@unibo.it}
\urladdr{http://www.lupini.org/}
\author{Martino Lupini}
\address{Dipartimento di Matematica, Universit\`{a} di Bologna, Piazza di
Porta S. Donato, 5, 40126 Bologna,\ Italy}
\email{martino.lupini@unibo.it}
\urladdr{http://www.lupini.org/}
\thanks{The authors were partially supported by the Marsden Fund Fast-Start Grant VUW1816 and the
Rutherford Discovery Fellowship VUW2002 from the Royal Society of New
Zealand, the Starting Grant 101077154
\textquotedblleft Definable Algebraic Topology\textquotedblright\ from the
European Research Council, the Gruppo
Nazionale per le Strutture Algebriche, Geometriche e le loro Applicazioni
(GNSAGA) of the Istituto Nazionale di Alta Matematica (INDAM), and the
University of Bologna.}
\subjclass[2000]{Primary 55P99, 55N07; Secondary 54H05, 55N10, 18G10, 18G35}
\keywords{Polish group, Steenrod property, group with a Polish cover, homogeneous space with a Polish cover, phantom subspaces, phantom length, abelian category, quasi-abelian category, left heart, homology, cohomology, Borel complexity, homotopy, solenoid, Steenrod homology, \v{C}ech cohomology, Steenrod duality, Alexander duality, phantom maps}
\date{\today }

\begin{abstract}
In this paper we isolate a general Universal Coefficient Theorem in the
context of abelian categories. We then apply it in the left heart of the
quasi-abelian category of abelian Polish groups to obtain short exact
sequences relating Steenrod homology and \v{C}ech cohomology, regarded as
functors to such a category. These are then used to (1) intrinsically
characterize the phantom subgroups of \v{C}ech cohomology via a recursive
purely algebraic formula, which can be see as a higher order generalization
of the Milnor exact sequence; (2) describe phantom subgroups of \v{C}ech
cohomology in terms of higher order phantom maps, in the context of the
homotopical description of \v{C}ech cohomology; (3) classify up to homotopy
(phantom) maps on higher order versions of solenoid complements, and measure
from the viewpoint of Borel complexity theory the complexity of such a
classification problem.
\end{abstract}

\maketitle

\section{Introduction}

In this paper, we isolate an abstract Universal Coefficient Theorem in the
context of abelian categories, subsuming the original Eilenberg--Steenrod
Universal Coefficient Theorems from \cite{eilenberg_group_1942} and their
algebraic version for complexes of modules \cite[Section 3.6]%
{weibel_introduction_1995}. We then apply it in the case of the left heart $%
\mathrm{LH}\left( \mathbf{PAb}\right) $ of the category $\mathbf{PAb}$ of
Polish abelian groups. The latter is a quasi-abelian category as defined in 
\cite{schneiders_quasi-abelian_1999}. Its left heart is an abelian category
uniquely characterized by a suitable universal property. While generally
constructed as a full subcategory of the derived category, an explicit
presentation of $\mathrm{LH}\left( \mathbf{PAb}\right) $ is obtained in \cite%
{lupini_looking_2024} in terms of abelian \emph{groups with a Polish cover }%
and \emph{Borel-definable homomorphism}.

In the context of this presentation, a canonical chain of subobjects of an
object of $\mathrm{LH}\left( \mathbf{PAb}\right) $, called its \emph{phantom
subgroups }(or Solecki subgroups), are introduced in \cite%
{lupini_looking_2024,lupini_complexity_2025} building on previous work of
Solecki \cite{solecki_equivalence_1995} and Farah--Solecki \cite%
{farah_borel_2006}.

Applying the general Universal Coefficient Theorem within $\mathrm{LH}\left( 
\mathbf{PAb}\right) $, we obtain exact sequences expressing Steenrod
homology and \v{C}ech cohomology, regarded as functors to $\mathrm{LH}\left( 
\mathbf{PAb}\right) $, one in terms of the other. Such exact sequences are
then applied to: (1) intrinsically describe the \emph{phantom subgroups }of 
\v{C}ech cohomology via a recursive purely algebraic formula; (2) refine the
homotopical description of cohomology by establishing a correspondence
between\emph{\ phantom subgroups} of \v{C}ech cohomology and \emph{higher
order phantom maps}; (3) classify up to homotopy maps from \emph{higher
order solenoid complements }to the $2$-dimensional sphere, and measure their
complexity from the viewpoint of Borel complexity theory.

Concerning these results: (1) can be seen as a higher order analogue of the
famous Milnor exact sequence \cite{milnor_axiomatic_1962,ferry_remarks_1995}%
, which corresponds to the particular case of the phantom subgroup of index
zero; (2) generalizes the well-known correspondence between the phantom
subgroup of index zero and phantom maps \cite[Section 7]%
{bergfalk_definable_2024-1}; (3) subsumes as a particular case the
Borsuk--Eilenberg Classification of maps on the $p$-adic solenoid complement 
\cite[Section 8]{bergfalk_definable_2024-1}, originally considered by Borsuk
and Eilenberg in 1936 \cite{borsuk_uber_1936}; see also \cite[Section 2.2]%
{kromer_tool_2007} and \cite{eilenberg_karol_1993,weibel_history_1999}. In
particular, we show that by replacing solenoids with their higher order
analogue one obtain classification problems that can have arbitrarily high
complexity from the perspective of Borel complexity theory \cite%
{gao_invariant_2009}.

This paper is divided into 8 sections, including this introduction. In
Section \ref{Section:homological-algebra} we briefly recall the necessary
background from category theory and homological algebra, including the
general UCT in the context of countably complete abelian categories. Section %
\ref{Section:topology} presents the construction of Steenrod homology and 
\v{C}ech cohomology with coefficients in an abelian category $\mathcal{M}$.
Section \ref{Section:complexity} recalls some notions from Borel complexity
theory, and the description of the left heart of abelian Polish groups in
terms of groups with a Polish cover. Section \ref{Section:phantom-maps}
presents phantom map, and their higher order version, and establishes some
of its fundamental properties. Section \ref{Section:Milnor} recalls the
Milnor Exact Sequence, and establishes its higher order analogue.\ Section %
\ref{Section:phantom} isolates consequences of the Universal Coefficient
Theorem and Milnor Exact Sequence for Steenrod Homology and \v{C}ech
cohomology. Finally, Section \ref{Section:higher} introduces higher order
analogues of solenoids and applies our main results to the homotopy
classification problem for maps on solenoid complements.

\subsubsection*{Acknowledgments}

We are grateful to Shaun Allison, Jeffrey Bergfalk, Oliver Braunling,
Alessandro Codenotti, Ivan Di Liberti, Nicholas Meadows, Marco Moraschini,
Aristotelis Panagiotopoulos, George Raptis, Filippo Sarti, and Joseph
Zielinski for many stimulating conversations and useful suggestions, and to
the President of Logic in\ Bologna for his unwavering commitment to
providing pastoral care to Members of the University Community.

\section{Homological algebra background\label{Section:homological-algebra}}

In this sections we recall notions from category theory and abstract
homological algebra to be used in the rest of the paper, as can be found in 
\cite%
{kashiwara_categories_2006,gelfand_methods_2003,weibel_introduction_1995,schneiders_quasi-abelian_1999,cartan_homological_1999,mac_lane_categories_1998,mac_lane_homology_1995}%
. In particular, we recall the notions of quasi-abelian and abelian
category, its corresponding derived category and left heart, as well as a
construction of derived functors in this context. The section ends with a
general Universal Coefficient Theorem in countably complete abelian
categories.

\subsection{The category of towers}

We begin with recalling the notion of \emph{tower} over a category, and the
category of towers. For a category $\mathcal{C}$, let $\mathcal{C}^{\vee }$
be the category of functors $\mathcal{C}^{\mathrm{op}}\rightarrow \mathbf{Set%
}^{\mathrm{op}}$, and consider the functor $\mathrm{k}_{\mathcal{C}}:%
\mathcal{C}\rightarrow \mathcal{C}^{\vee }$ defined by $X\mapsto \mathrm{Hom}%
\left( X,-\right) $ \cite[Definition 1.4.2]{kashiwara_categories_2006}. A
version of the Yoneda Lemma asserts that $\mathrm{k}_{\mathcal{C}}$ is a
fully faithful functor, and in fact for every object $B$ of $\mathcal{C}%
^{\vee }$ and $X$ of $\mathcal{C}$ there is a bijection $\mathrm{Hom}_{%
\mathcal{C}^{\vee }}\left( B,\mathrm{k}_{\mathcal{C}}\left( X\right) \right)
\simeq B\left( X\right) $ which is natural in $B$ and $X$ \cite[Proposition
1.4.3]{kashiwara_categories_2006}. This allows one to identify $\mathcal{C}$
as a full subcategory of $\mathcal{C}^{\vee }$. The category $\mathcal{C}%
^{\vee }$ has small limits and colimits which are computed \textquotedblleft
levelwise\textquotedblright , namely for a functor $F:I\rightarrow \mathcal{C%
}^{\vee }$ we have that 
\begin{equation*}
\left( \mathrm{\mathrm{lim}}_{i}F\left( i\right) \right) \left( X\right)
\cong \mathrm{\mathrm{lim}}_{i}\left( F\left( i\right) \left( X\right)
\right) \text{ and }\left( \mathrm{col\mathrm{im}}_{i}F\left( i\right)
\right) \left( X\right) \cong \mathrm{col\mathrm{im}}_{i}\left( F\left(
i\right) \left( X\right) \right) \text{;}
\end{equation*}%
see \cite[Corollary 2.4.3]{kashiwara_categories_2006}. When $\mathcal{C}$
admits small colimits, the inclusion functor $\mathcal{C}\rightarrow 
\mathcal{C}^{\vee }$ is cocontinuous, namely preserves small colimits.
However, in general $\mathrm{k}_{\mathcal{C}}$ is not continuous (i.e., does
not preserve small limits). Thus for a functor $F:I\rightarrow \mathcal{C}$
one sets $\mathrm{\mathrm{lim}}^{\vee }F$ to be the limit in $\mathcal{C}%
^{\vee }$ of the functor $F:I\rightarrow \mathcal{C}^{\vee }$. This is
called the \emph{pro-limit }of $F$; see \cite[Notation 2.6.1]%
{kashiwara_categories_2006}.

A \emph{tower }over $\mathcal{C}$ is an object of $\mathcal{C}^{\vee }$ that
is isomorphic to a pro-limit of a functor $F:\omega ^{\mathrm{op}%
}\rightarrow \mathcal{C}$ where $\omega $ is the linear order of natural
numbers \cite[Definition 6.1.1]{kashiwara_categories_2006}. The \emph{%
category of towers }\textrm{Pro}$_{\omega }\left( \mathcal{C}\right) $ is
the full subcategory of $\mathcal{C}^{\vee }$ spanned by towers. One can
explicitly describe morphisms in $\mathrm{Pro}_{\omega }\left( \mathcal{C}%
\right) $ as in \cite[Definition 2.1.2]{edwards_cech_1976}. If $\mathcal{C}$
is finitely complete and finitely cocomplete, then $\mathrm{Pro}_{\omega
}\left( \mathcal{C}\right) $ is countably complete and finitely cocomplete 
\cite[Section 2]{blanc_colimits_1996}. Dually, the category $\mathrm{Ind}%
_{\omega }\left( \mathcal{C}\right) $, whose objects are called inductive
sequences over $\mathcal{C}$, is defined by setting $\mathrm{Ind}_{\omega
}\left( \mathcal{C}\right) :=\mathrm{Pro}_{\omega }\left( \mathcal{C}^{%
\mathrm{op}}\right) ^{\mathrm{op}}$. When $\mathcal{C}$ is countably
continuous, the inclusion $j_{\mathcal{C}}:\mathcal{C}\rightarrow \mathrm{Pro%
}_{\omega }\left( \mathcal{C}\right) $ is not countably continuous in
general. Taking the limit of a tower defines a functor $\mathrm{\mathrm{lim}}%
:\mathrm{Pro}_{\omega }\left( \mathcal{C}\right) \rightarrow \mathcal{C}$
that is the right adjoint of $j_{C}$ and such that $\mathrm{\mathrm{lim}}%
\circ j_{C}$ is isomorphic to the identity.

\subsection{Abelian and quasi-abelian categories}

We now recall some notions concerning \emph{abelian} categories and, more
generally, \emph{additive} and \emph{quasi-abelian} categories as can be
found in \cite{mac_lane_categories_1998,schneiders_quasi-abelian_1999}. An%
\emph{\ }$\mathbf{Ab}$\textrm{-}\emph{category} is a category $\mathcal{A}$
enriched over the category $\mathbf{Ab}$ of abelian groups. Thus, for each
pair of objects $X,Y$ of $\mathcal{A}$, $\mathrm{Hom}_{\mathcal{A}}\left(
X,Y\right) $ is an abelian group, in such a way that composition of arrows 
\textrm{Hom}$_{\mathcal{A}}\left( Y,Z\right) \times \mathrm{Hom}_{\mathcal{A}%
}\left( X,Y\right) \rightarrow \mathrm{Hom}_{\mathcal{A}}\left( X,Z\right) $
is bilinear for each triple of objects $X,Y,Z$ in $\mathcal{A}$. In an $%
\mathbf{Ab}$\textrm{-}category, an object $X$ is terminal if and only if it
is initial if and only if $\mathrm{Hom}_{\mathcal{A}}\left( X,X\right) $ is
the trivial group, in which case $X$ is called a \emph{null object }for $%
\mathcal{A}$. Furthermore, binary products and binary coproducts coincide
and are called \emph{biproducts}. An \emph{additive category} is an $\mathbf{%
Ab}$-category that has biproducts and a null object. A \emph{quasi-abelian
category }is an additive category that has kernels and cokernels (and hence
all finite limits and finite colimits), and such that the class of kernels
is closed under push-out along arbitrary morphisms and the class of
cokernels is closed under pull-back along arbitrary morphisms. An \emph{%
abelian category} is an additive category where every monic is a kernel and
every epic is a cokernel.

In a quasi-abelian category $\mathcal{A}$, one can define the image $\mathrm{%
im}\left( f\right) $ of an arrow to be $\mathrm{ker}\left( \mathrm{coker}%
\left( f\right) \right) $, and its coimage to be $\mathrm{coker}\left( 
\mathrm{ker}\left( f\right) \right) $. Then $f$ induces a canonical map $%
\overline{f}:\mathrm{coim}\left( f\right) \rightarrow \mathrm{im}\left(
f\right) $ which is both monic and epic. The arrow is called \emph{strict }%
if $\overline{f}$ is an isomorphism. Then we have that an arrow is monic and
strict if and only if it is a kernel, and it is epic and strict if and only
if it is a cokernel. An abelian category is precisely a quasi-abelian
category where every arrow is strict.

A functor $F:\mathcal{A}\rightarrow \mathcal{B}$ between quasi-abelian
categories is:

\begin{itemize}
\item \emph{additive }if it induces group homomorphism $\mathrm{Hom}_{%
\mathcal{A}}\left( X,Y\right) \rightarrow \mathrm{Hom}_{\mathcal{B}}\left(
FX,FY\right) $ for each pair of objects $X,Y$ in $\mathcal{A}$ or,
equivalently, it preserves biproducts;

\item \emph{left exact }if it preserves kernels of strict arrows;

\item \emph{right exact }if it preserves cokernels of strict arrows;

\item \emph{exact }if it is both left and right exact.
\end{itemize}

\subsection{Complexes and the derived category}

An important construction that one can perform starting from an abelian
category is to take its \emph{derived category }\cite[Chapter III]%
{gelfand_methods_2003}, which is the natural setting for the study of \emph{%
derived functors}. More generally, such a construction can be performed
starting from a quasi-abelian category $\mathcal{A}$, as explained in \cite%
{schneiders_quasi-abelian_1999}. We now recall the main ingredients in this
construction.

A \emph{complex }over $\mathcal{A}$ is a sequence%
\begin{equation*}
\cdots \longrightarrow A_{2}\overset{\partial _{2}}{\longrightarrow }A_{1}%
\overset{\partial _{1}}{\longrightarrow }A_{0}\overset{\partial _{0}}{%
\longrightarrow }A_{-1}\longrightarrow \cdots
\end{equation*}%
of objects $A_{n}$ of $\mathcal{A}$ and morphisms $\partial
_{n}:A_{n}\rightarrow A_{n-1}$ of $\mathcal{A}$ for $n\in \mathbb{Z}$
satisfying $\partial _{n-1}\partial _{n}=0$. Given such a complex $A$, one
sets $A^{n}:=A_{-n}$ and $\delta ^{n}:=\partial _{-n}:A^{n}\rightarrow
A^{n+1}$. A chain complex is a complex with $A_{n}=0$ for $n<0$ and a
cochain complex is a complex with $A^{n}=0$ for $n<0$. A map $f:A\rightarrow
B$ between complexes is a sequence $\left( f_{n}\right) _{n\in \omega }$ of
morphisms in $\mathcal{A}$ satisfying $\partial f_{n}=f_{n-1}\partial $ for $%
n\in \mathbb{Z}$. A map $f$ is\emph{\ null-homotopic} if there exist
morphisms $s_{n}:A_{n}\rightarrow B_{n+1}$ in $\mathcal{A}$ such that $%
\partial s_{n}+s_{n-1}\partial =f_{n}$ for every $n\in \mathbb{Z}$, while
two maps $f$ and $g$ are homotopic if $f-g$ is null-homotopic. Complexes
over $\mathcal{A}$ form a category $\mathbf{Ch}\left( \mathcal{A}\right) $
that has maps of complexes as morphisms. The category $\mathbf{K}\left( 
\mathcal{A}\right) $ has complexes as objects and\emph{\ homotopy classes }%
of maps as morphisms. The latter category is not quasi-abelian in general,
even when $\mathcal{A}$ is abelian.\ It has, however, the structure of \emph{%
triangulated category }\cite[Chapter 10]{kashiwara_categories_2006}.

A complex $A$ in $\mathbf{K}\left( \mathcal{A}\right) $ is exact if\emph{\ }$%
\partial _{n}^{A}:A_{n}\rightarrow A_{n-1}$ is strict and $\mathrm{k\mathrm{%
er}}\left( \partial _{n}^{A}\right) =\mathrm{im}\left( \partial
_{n+1}^{A}\right) $ for every $n\in \mathbb{Z}$. The full subcategory $%
\mathbf{N}\left( \mathcal{A}\right) $ of $\mathbf{K}\left( \mathcal{A}%
\right) $ spanned by exact complexes is a null system \cite[Definition 10.2.2%
]{kashiwara_categories_2006}. A map $f:A\rightarrow B$ is a \emph{%
quasi-isomorphism} if its mapping cone $\mathrm{cone}\left( f\right) $ is an
exact complex \cite[Section 1.5]{weibel_introduction_1995}. The subcategory $%
S_{\mathbf{N}\left( \mathcal{A}\right) }$ of $\mathbf{K}\left( \mathcal{A}%
\right) $ whose arrows are the quasi-isomorphisms is a saturated
multiplicative system in $\mathbf{K}\left( \mathcal{A}\right) $, and the
corresponding quotient $\mathbf{D}\left( \mathcal{A}\right) =\mathbf{K}%
\left( \mathcal{A}\right) [S_{\mathbf{N}\left( \mathcal{A}\right) }^{-1}]$
is a triangulated category, called the\emph{\ derived category} of $\mathcal{%
A}$. Identifying an object $A$ of $\mathcal{A}$ with the complex with $%
A^{0}=A$ and $A^{n}=0$ for $n\neq 0$ allows one to regard $\mathcal{A}$ as a
full subcategory of $\mathbf{D}\left( \mathcal{A}\right) $.

\subsection{The left heart of a quasi-abelian category}

Suppose that $\mathcal{A}$ is a quasi-abelian category. Then the derived
category $\mathbf{D}\left( \mathcal{A}\right) $ is endowed with a canonical
left\emph{\ truncation structure} (t-structure) \cite{kashiwara_sheaves_1994}%
. The heart of this t-structure on $\mathbf{D}\left( \mathcal{A}\right) $ is
called, for brevity, the \emph{left heart }of $\mathcal{A}$ and denoted by $%
\mathrm{LH}\left( \mathcal{A}\right) $. This is a full subcategory of $%
\mathbf{D}\left( \mathcal{A}\right) $ containing $\mathcal{A}$, such that
the inclusion $\mathcal{A}\rightarrow \mathrm{LH}\left( \mathcal{A}\right) $
is exact and finitely continuous. Furthermore, $\mathrm{LH}\left( \mathcal{A}%
\right) $ is an abelian category, and it is characterized by the following
universal property: for every exact and finitely continuous functor $F:%
\mathcal{A}\rightarrow \mathcal{M}$, where $\mathcal{M}$ is an abelian
category, there exists an essentially unique exact and finitely continuous
functor $\mathrm{LH}\left( \mathcal{A}\right) \rightarrow \mathcal{M}$ whose
restriction to $\mathcal{A}$ is isomorphic to $F$. The inclusion $\mathcal{A}%
\rightarrow \mathrm{LH}\left( \mathcal{A}\right) $ extends to an equivalence
of categories $\mathbf{D}\left( \mathcal{A}\right) \rightarrow \mathbf{D}%
\left( \mathrm{LH}\left( \mathcal{A}\right) \right) $. When $\mathcal{A}$ is
abelian, one has that $\mathcal{A}=\mathrm{LH}\left( \mathcal{A}\right) $.

The elements of $\mathrm{LH}\left( \mathcal{A}\right) $ are precisely the
complexes over $A$ that are isomorphic in $\mathbf{D}\left( \mathcal{A}%
\right) $ to a complex of the form%
\begin{equation*}
\cdots \longrightarrow 0\longrightarrow A_{1}\overset{\partial _{1}^{A}}{%
\longrightarrow }A_{0}\longrightarrow 0\longrightarrow \cdots 
\end{equation*}%
where $\partial _{1}^{A}$ is a monic in $\mathcal{A}$. The homology $\mathrm{%
H}_{\bullet }\left( A\right) $ of a complex $A$ over $\mathcal{A}$ is a
graded object of $\mathrm{LH}\left( \mathcal{A}\right) $, namely an element
of the functor category $\mathrm{LH}\left( \mathcal{A}\right) ^{\mathbb{Z}}$
where $\mathbb{Z}$ is regarded as a category with only identity arrows. This
is defined by setting, for $n\in \mathbb{Z}$, \textrm{H}$_{n}\left( A\right) 
$ to be the element of $\mathrm{LH}\left( \mathcal{A}\right) $ represented
by complex%
\begin{equation*}
0\longrightarrow \mathrm{coim}\left( \partial _{n}^{A}\right) \overset{f}{%
\longrightarrow }\mathrm{ker}\left( \partial _{n-1}^{A}\right)
\longrightarrow 0
\end{equation*}%
where $f$ is the canonical morphism induced by $\partial
_{n}^{A}:A_{n}\rightarrow A_{n-1}$ by the universal properties of the kernel
and the cokernel. This defines a functor $\mathrm{H}_{\bullet }:\mathbf{D}%
\left( \mathcal{A}\right) \rightarrow \mathrm{LH}\left( \mathcal{A}\right) ^{%
\mathbb{Z}}$. Dually, one defines the cohomology $\mathrm{H}^{n}\left(
A\right) $ to be $\mathrm{H}_{-n}\left( A\right) $. One has that a map $%
f:A\rightarrow B$ between complexes over $\mathcal{A}$ is a \emph{%
quasi-isomorphism} if and only if it induces an \emph{isomorphism} at the
level of homology. Two morphisms of complexes are \emph{cohomologous }if
they induce the same map in cohomology.

\subsection{Derived functors\label{Subsection:derived}}

Suppose that $\mathcal{A}$ is an \emph{abelian} category. An object $E$ of $%
\mathcal{A}$ is \emph{injective }if the functor $\mathrm{Hom}_{\mathcal{A}%
}\left( -,E\right) $ is exact \cite[Lemma 2.2.3]{weibel_introduction_1995}.
The category $\mathcal{A}$ \emph{has enough injectives }if for every object $%
A$ of $\mathcal{A}$ there exists a monic $A\rightarrow E$ where $E$ is
injective. Suppose that $\mathcal{A}$ has enough injectives, and let $%
\mathcal{I}$ be the (additive) full subcategory of $\mathcal{A}$ spanned by
the injective objects. Let $\mathbf{K}^{+}\left( \mathcal{A}\right) $
(respectively, $\mathbf{K}^{+}\left( \mathcal{I}\right) $) be the
triangulated categories of \emph{bounded below} complexes over $\mathcal{A}$
(respectively, over $\mathcal{I}$) and homotopy classes of maps \cite[%
Corollary 10.2.5]{weibel_introduction_1995}. Let also $\mathbf{D}^{+}\left( 
\mathcal{A}\right) \mathbf{\ }$be the triangulated category $\mathbf{K}%
^{+}\left( \mathcal{A}\right) /\mathbf{N}^{+}\left( \mathcal{A}\right) $
where $\mathbf{N}^{+}\left( \mathcal{A}\right) $ is the null system of
bounded below exact complexes over $\mathcal{A}$, and let $Q_{\mathbf{D}%
^{+}\left( \mathcal{A}\right) }:\mathbf{K}^{+}\left( \mathcal{A}\right)
\rightarrow \mathbf{D}^{+}\left( \mathcal{A}\right) $ be the canonical
quotient functor, which is finitely continuous and finitely cocontinuous 
\cite[Corollary 10.4.3]{weibel_introduction_1995}. The inclusion $\mathcal{I}%
\rightarrow \mathcal{A}$ induces a triangulated functor $\mathbf{K}%
^{+}\left( \mathcal{I}\right) \rightarrow \mathbf{K}^{+}\left( \mathcal{A}%
\right) $ that composed with the canonical functor $\mathbf{K}^{+}\left( 
\mathcal{A}\right) \rightarrow \mathbf{D}^{+}\left( \mathcal{A}\right) $
yields an equivalence of triangulated categories \cite[Definition 10.1.9]%
{kashiwara_categories_2006}; see \cite[Theorem 10.4.8]%
{weibel_introduction_1995}.

Suppose that $\mathcal{M}$ is an abelian category, and $F:\mathcal{I}%
\rightarrow \mathcal{M}$ is an additive functor. The \emph{right derived
functor }of $F$ is the functor\textbf{\ }$\mathbf{R}F:\mathbf{D}^{+}\left( 
\mathcal{A}\right) \rightarrow \mathbf{D}^{+}\left( \mathcal{M}\right) $
obtained as a composition%
\begin{equation*}
\mathbf{D}^{+}\left( \mathcal{A}\right) \rightarrow \mathbf{K}^{+}\left( 
\mathcal{I}\right) \rightarrow \mathbf{K}^{+}\left( \mathcal{M}\right)
\rightarrow \mathbf{D}^{+}\left( \mathcal{M}\right)
\end{equation*}%
where $\mathbf{D}^{+}\left( \mathcal{A}\right) \rightarrow \mathbf{K}%
^{+}\left( \mathcal{I}\right) $ is the inverse of the triangulated
equivalence described above, $\mathbf{K}^{+}\left( \mathcal{I}\right)
\rightarrow \mathbf{K}^{+}\left( \mathcal{M}\right) $ is the triangulated
functor induced by $F$, and $\mathbf{K}^{+}\left( \mathcal{M}\right)
\rightarrow \mathbf{D}^{+}\left( \mathcal{M}\right) $ is the quotient
functor; see \cite[Proposiotion 1.3.5]{schneiders_quasi-abelian_1999}. Such
a functor is characterized by a suitable universal property; see \cite[%
Definition 13.3.1]{kashiwara_categories_2006} and \cite[Definition 1.3.1]%
{schneiders_quasi-abelian_1999}. In particular, there exists an essentially
unique finitely continuous functor $\mathcal{A}\rightarrow \mathcal{M}$,
still denoted by $F$, whose restriction to $\mathcal{I}$ is isomorphic to $F$%
, which can be simply defined as $\left( \mathrm{H}^{0}\circ \mathbf{R}%
F\right) |_{\mathcal{A}}$. One also has that such an extension is \emph{%
countably }continuous when the original functor $F:\mathcal{I}\rightarrow 
\mathcal{M}$ is \emph{countably} continuous. Defining $R^{n}F:=\left( 
\mathrm{H}^{n}\circ \mathbf{R}F\right) |_{\mathcal{A}}$ for $n\geq 0$ one
obtains a \emph{universal cohomological }$\delta $-functor with $R^{0}F|_{%
\mathcal{A}}\cong F$ \cite[Definition 2.1.1, Definition 2.1.4, Theorem 2.4.7]%
{weibel_introduction_1995}. In particular, for every short-exact sequence%
\begin{equation*}
0\rightarrow A\rightarrow B\rightarrow C\rightarrow 0
\end{equation*}%
in $\mathcal{A}$, there is a corresponding long exact sequence%
\begin{equation*}
0\rightarrow FA\rightarrow FB\rightarrow FC\rightarrow R^{1}F\left( A\right)
\rightarrow R^{1}F\left( B\right) \rightarrow R^{1}F\left( C\right)
\rightarrow R^{2}F\left( A\right) \rightarrow \cdots
\end{equation*}%
in $\mathcal{M}$, where the morphisms%
\begin{equation*}
R^{n}F\left( A\right) \rightarrow R^{n+1}F\left( C\right)
\end{equation*}%
are natural transformations.

\subsection{The derived functor of \textrm{Hom}\label{Subsection:derive-hom}}

Consider the category $\mathbf{Ab}_{\aleph _{0}}$ of countable discrete
abelian groups and let $\mathcal{P}$ be the full subcategory of $\mathbf{Ab}%
_{\aleph _{0}}$ spanned by the countable free abelian groups. Notice that
the objects of $\mathcal{P}$ are precisely the projective objects in $%
\mathbf{Ab}_{\aleph _{0}}$ (namely, the injective objects in $\mathbf{Ab}%
_{\aleph _{0}}^{\mathrm{op}}$). Suppose that $\mathcal{M}$ is a countably
complete abelian category, and $X$ is an object in $\mathcal{M}$.

We use the notation of generalized \textquotedblleft
elements\textquotedblright\ of objects of $\mathcal{M}$ as in \cite[Section
VIII.4]{mac_lane_categories_1998} to describe morphisms in $\mathcal{M}$.
For an element $x$ of a product $\prod_{i\in I}A_{i}$ of objects of $%
\mathcal{M}$, we let $x_{i}$ for $i\in I$ be the images of $x$ under the
canonical projections. If $\mathbb{Z}^{\left( I\right) }$ is the free
abelian group on a countable set $I$ with free basis $\left( e_{i}\right)
_{i\in I}$, we set%
\begin{equation*}
\mathrm{Hom}(\mathbb{Z}^{\left( I\right) },X):=\prod\nolimits_{I}X\text{,}
\end{equation*}%
where the product is computed in $\mathcal{M}$. Suppose that $I$ and $J$ are
countable sets, and let $\varphi :\mathbb{Z}^{\left( I\right) }\rightarrow 
\mathbb{Z}^{\left( J\right) }$ be the homomorphism given by%
\begin{equation*}
e_{i}\mapsto \sum_{j\in J}a_{ij}e_{j}
\end{equation*}%
for $a_{ij}\in \mathbb{Z}$. We define%
\begin{equation*}
\mathrm{Hom}\left( \varphi ,X\right) :\prod\nolimits_{J}X\rightarrow
\prod\nolimits_{I}X
\end{equation*}%
to be the morphism in $\mathcal{M}$ such that, for every element $x$ of $%
\prod\nolimits_{J}X$ and $i\in I$,%
\begin{equation*}
\mathrm{Hom}\left( \varphi ,X\right) \left( x\right) _{i}=\sum_{j\in
J}a_{ij}x_{j}\text{.}
\end{equation*}%
Notice that this is well-defined since, for every $i\in I$, $\left\{ j\in
J:a_{ij}\neq 0\right\} $ is finite.

This defines an additive\emph{\ }functor $\mathrm{Hom}\left( -,X\right) :%
\mathcal{P}^{\mathrm{op}}\rightarrow \mathcal{M}$. By the discussion in
Section \ref{Subsection:derived} applied to the case when $\mathcal{A}=%
\mathbf{Ab}_{\aleph _{0}}^{\mathrm{op}}$ and $F=\mathrm{Hom}\left(
-,X\right) $, we have that $\mathrm{Hom}\left( -,X\right) $ has an
essentially unique extension to a finitely continuous functor $\mathrm{Hom}%
\left( -,X\right) :\mathbf{Ab}_{\aleph _{0}}^{\mathrm{op}}\rightarrow 
\mathcal{M}$, which is in fact countably continuous. Furthermore, we have
that $\mathrm{Hom}\left( -,X\right) $ has a right derived functor 
\begin{equation*}
\mathbf{Ext}\left( -,X\right) :\mathbf{D}^{+}\left( \mathbf{Ab}_{\aleph
_{0}}^{\mathrm{op}}\right) \rightarrow \mathbf{D}^{-}\left( \mathcal{M}%
\right) \text{.}
\end{equation*}%
Setting%
\begin{equation*}
\mathrm{Ext}\left( -,X\right) :=\left( \mathrm{H}^{1}\circ \mathbf{Ext}_{_{%
\mathcal{M}}}\left( -,X\right) \right) |_{\mathbf{Ab}^{\mathrm{op}}}
\end{equation*}%
one obtains a functor $\mathrm{Ext}\left( -,X\right) :\mathbf{Ab}_{\aleph
_{0}}^{\mathrm{op}}\rightarrow \mathcal{M}$. Every short-exact sequence%
\begin{equation*}
0\rightarrow A\rightarrow B\rightarrow C\rightarrow 0
\end{equation*}%
of abelian groups gives rise to an exact sequence%
\begin{equation*}
0\rightarrow \mathrm{Hom}\left( C,X\right) \rightarrow \mathrm{Hom}\left(
B,X\right) \rightarrow \mathrm{Hom}\left( A,X\right) \rightarrow \mathrm{Ext}%
\left( C,X\right) \rightarrow \mathrm{Ext}\left( B,X\right) \rightarrow 
\mathrm{Ext}\left( A,X\right) \rightarrow 0\text{.}
\end{equation*}

Suppose that $A$ is a complex of countable free abelian groups. This gives
rise to a complex $\mathrm{Hom}\left( A,X\right) $ in $\mathcal{M}$,
obtained by setting 
\begin{equation*}
\mathrm{Hom}\left( A,X\right) _{n}:=\mathrm{Hom}\left( A^{n},X\right)
\end{equation*}%
for every $n\in \mathbb{Z}$. The usual proof of the Universal Coefficient
Theorem \cite[Theorem 3.6.5]{weibel_introduction_1995} for modules gives the
following.

\begin{theorem}
\label{Theorem:UCT}Suppose that $A$ is a complex of\emph{\ countable} \emph{%
free} abelian groups. Let also $\mathcal{M}$ be a countably complete abelian
category, and $X$ be an object of $\mathcal{M}$. Then we have a natural
short exact sequence of graded objects in $\mathcal{M}$%
\begin{equation*}
0\rightarrow \mathrm{Ext}\left( \mathrm{H}^{\bullet +1}(A),X\right)
\rightarrow \mathrm{H}_{\bullet }\left( \mathrm{Hom}\left( A,X\right)
\right) \rightarrow \mathrm{Hom}\left( \mathrm{H}^{\bullet }\left( A\right)
,X\right) \rightarrow 0
\end{equation*}%
that splits (unnaturally). This is called the \emph{Universal Coefficient
Sequence }for $A$ and $X$.
\end{theorem}

\begin{corollary}
\label{Corollary:UCT}The assignment $A\mapsto \mathrm{Hom}\left( A,X\right) $
induces a functor $\mathbf{D}\left( \mathbf{Ab}_{\aleph _{0}}\right) ^{%
\mathrm{op}}\rightarrow \mathbf{D}\left( \mathcal{M}\right) $.
\end{corollary}

\begin{proof}
It suffices to notice that a quasi-isomorphism in $\mathbf{K}\left( \mathbf{%
Ab}\right) $ induces a quasi-isomorphism in $\mathbf{K}\left( \mathcal{M}%
\right) $ by Theorem \ref{Theorem:UCT} and the Five Lemma \cite[Lemma 8.3.13]%
{kashiwara_categories_2006}.
\end{proof}

\begin{remark}
\label{Remark:pairs}Suppose that 
\begin{equation*}
0\rightarrow A\rightarrow B\rightarrow C\rightarrow 0
\end{equation*}%
is a short exact sequence of complexes of countable abelian groups. We say
that it is \emph{locally split }if%
\begin{equation*}
0\rightarrow A_{n}\rightarrow B_{n}\rightarrow C_{n}\rightarrow 0
\end{equation*}%
is a split short exact sequence of countable abelian groups for every $n\in 
\mathbb{Z}$. In this case, we have a corresponding short exact sequence%
\begin{equation*}
0\rightarrow \mathrm{Hom}\left( C,X\right) \rightarrow \mathrm{Hom}\left(
B,X\right) \rightarrow \mathrm{Hom}\left( A,X\right) \rightarrow 0\text{.}
\end{equation*}%
This induces a long exact sequence in homology%
\begin{equation*}
\cdots \rightarrow \mathrm{H}_{n}\left( \mathrm{Hom}\left( C,X\right)
\right) \rightarrow \mathrm{H}_{n}\left( \mathrm{Hom}\left( B,X\right)
\right) \rightarrow \mathrm{H}_{n}\left( \mathrm{Hom}\left( A,X\right)
\right) \overset{d_{n}}{\rightarrow }\mathrm{H}_{n-1}\left( \mathrm{Hom}%
\left( C,X\right) \right) \rightarrow \cdots
\end{equation*}%
see \cite[Theorem 1.3.1]{weibel_introduction_1995}. Likewise, we have a long
exact sequence in cohomology%
\begin{equation*}
\cdots \rightarrow \mathrm{H}^{n}\left( C\right) \rightarrow \mathrm{H}%
^{n}\left( B\right) \rightarrow \mathrm{H}^{n}\left( C\right) \rightarrow 
\mathrm{H}^{n}\left( A\right) \overset{d^{n}}{\rightarrow }\mathrm{H}%
^{n+1}\left( C\right) \rightarrow \cdots
\end{equation*}%
One can see that the connecting morphisms $d_{n}$ and $d^{n}$ define a \emph{%
natural transformation} from the Universal Coefficient Sequence for $A$ and $%
X$ to the Universal Coefficient Sequence for $C$ and $X$.
\end{remark}

\subsection{Homotopy colimits}

We now recall the construction of \emph{homotopy limits} and\emph{\ homotopy
colimits }of complexes over an abelian category, to be used in the
definition of homology and cohomology groups of spaces; see \cite[Section 17]%
{mardesic_strong_2000}. Suppose that $\mathcal{B}$ is a countably cocomplete
abelian category. Let 
\begin{equation*}
\boldsymbol{A}=(A_{\left( n\right) },\eta _{\left( m,n\right) }:A_{\left(
n\right) }\rightarrow A_{\left( m\right) })
\end{equation*}%
be an inductive sequence of complexes over $\mathcal{B}$. One defines the 
\emph{homotopy colimit }of $\boldsymbol{A}$ to be the complex $\mathrm{hoco%
\mathrm{lim}}\boldsymbol{A}$ over $\mathcal{B}$ defined by setting%
\begin{equation*}
\left( \mathrm{hoco\mathrm{lim}}\boldsymbol{A}\right) ^{n}:=\coprod_{\ell
\in \omega }A_{\left( \ell \right) }^{n+1}\oplus \coprod_{\ell \in \omega
}A_{\left( \ell \right) }^{n}
\end{equation*}%
For $\ell \in \omega $, generalized elements $z$ of $A_{\left( \ell \right)
}^{n}$, and $w$ of $A_{\left( \ell \right) }^{n+1}$, we define the
generalized element $ze_{\ell }$ of $\left( \mathrm{hoco\mathrm{lim}}%
\boldsymbol{A}\right) ^{n}$ to be the image of $z$ under the canonical
morphism $A_{\left( \ell \right) }^{n}\rightarrow \left( \mathrm{hoco\mathrm{%
lim}}\boldsymbol{A}\right) ^{n}$, and the generalized element $we_{\ell
,\ell +1}$ of $\left( \mathrm{hoco\mathrm{lim}}\boldsymbol{A}\right) ^{n}$
to be the image of $w$ under the canonical morphism $A_{\left( \ell \right)
}^{n+1}\rightarrow \left( \mathrm{hoco\mathrm{lim}}\boldsymbol{A}\right)
^{n} $. We define the codifferential $\delta ^{n}:\left( \mathrm{hoco\mathrm{%
lim}}\boldsymbol{A}\right) ^{n}\rightarrow \left( \mathrm{hoco\mathrm{lim}}%
\boldsymbol{A}\right) ^{n+1}$ to be the morphism defined by setting, for $%
\ell \in \omega $, and generalized elements $z$ and $w$ of $A_{\left( \ell
\right) }^{n}$ and $A_{\left( \ell \right) }^{n+1}$, respectively,%
\begin{equation*}
\delta ^{n}\left( ze_{\ell }\right) =\left( \delta ^{n}z\right) e_{\ell }
\end{equation*}%
and%
\begin{equation*}
\delta ^{n}\left( we_{\ell ,\ell +1}\right) =\left( \delta ^{n+1}w\right)
e_{\ell ,\ell +1}+\left( -1\right) ^{n}(we_{\ell }-\eta _{\left( \ell
+1,\ell \right) }(w)e_{\ell +1})\text{;}
\end{equation*}%
see \cite[Section 2]{alonso_tarrio_localization_2000}.

Suppose now that $\boldsymbol{A}$ and $\boldsymbol{B}$ are inductive
sequences of complexes over $\mathcal{B}$, and $f=(t_{\ell },f_{\left( \ell
\right) })$ represents a morphism from $\boldsymbol{A}$ to $\boldsymbol{B}$
in $\mathrm{Ind}_{\omega }\left( \mathbf{Ch}\left( \mathcal{B}\right)
\right) $. This means that $\left( t_{\ell }\right) $ is an increasing
sequence in $\omega $ and $f_{\left( \ell \right) }:A_{\left( \ell \right)
}\rightarrow B_{\left( t_{\ell }\right) }$ is a map such that $\eta _{\left(
t_{\ell +1},t_{\ell }\right) }f_{\left( \ell \right) }=f_{\left( \ell
+1\right) }\eta _{\left( \ell +1,\ell \right) }$ for every $\ell \in \omega $%
. Then $f$ induces a map of complexes $f_{\left( \infty \right) }:=\mathrm{%
hoco\mathrm{lim}}f:\mathrm{hoco\mathrm{lim}}\boldsymbol{A}\rightarrow 
\mathrm{hoco\mathrm{lim}}\boldsymbol{B}$, defined in terms of generalized
elements by setting, for $n\in \mathbb{Z}$, $\ell \in \omega $, and
generalized elements $z$ and $w$ of $A_{\left( \ell \right) }^{n}$ and $%
A_{\left( \ell \right) }^{n+1}$, respectively, 
\begin{equation*}
f_{\left( \infty \right) }^{n}\left( ze_{\ell }\right) =f_{\left( \ell
\right) }^{n}\left( z\right) e_{t_{\ell }}
\end{equation*}%
and 
\begin{equation*}
f_{\left( \infty \right) }^{n}\left( we_{\ell ,\ell +1}\right) =f_{\left(
\ell \right) }^{n+1}(w)e_{t_{\ell },t_{\ell +1}}\text{.}
\end{equation*}%
This defines a functor from $\mathrm{hoco\mathrm{lim}}:\mathrm{Ind}_{\omega
}\left( \mathbf{Ch}\left( \mathcal{B}\right) \right) \rightarrow \mathbf{K}%
\left( \mathcal{B}\right) $.

Adopting the notation above, for every $\ell \in \omega $ we have a map%
\begin{equation*}
A_{\left( \ell \right) }\rightarrow \mathrm{hocolim}\boldsymbol{A}\text{, }%
z\mapsto ze_{\ell }\text{.}
\end{equation*}%
This induces a morphism%
\begin{equation*}
\mathrm{H}^{\bullet }(A_{\left( \ell \right) })\rightarrow \mathrm{H}%
^{\bullet }\left( \mathrm{hocolim}\boldsymbol{A}\right) \text{.}
\end{equation*}%
Such morphisms for $\ell \in \omega $ induce a morphism%
\begin{equation*}
\mathrm{co\mathrm{lim}}_{\ell }H^{\bullet }(A_{\left( \ell \right)
})\rightarrow H^{\bullet }\left( \mathrm{hocolim}\boldsymbol{A}\right) \text{%
,}
\end{equation*}%
which can be seen to be a natural isomorphism.

\begin{lemma}
\label{Lemma:cohomology-colimit}Suppose that $\boldsymbol{A}$ is an
inductive sequence of complexes over a countably cocomplete abelian category 
$\mathcal{B}$. Then $\mathrm{H}^{\bullet }\left( \mathrm{hoco\mathrm{lim}}%
\boldsymbol{A}\right) $ and $\mathrm{co\mathrm{lim}}_{\ell }{}\mathrm{H}%
^{\bullet }(A_{\left( \ell \right) })$ are naturally isomorphic.
\end{lemma}

The following lemmas provide sufficient conditions for homotopy colimits of
inductive sequences of complexes to be quasi-isomorphic.

\begin{lemma}
\label{Lemma:qi-1}Suppose that $\boldsymbol{A}$ and $\boldsymbol{B}$ are
inductive sequences of complexes over $\mathcal{B}$. Suppose that there
exist maps $f_{\left( \ell \right) }:B_{\left( \ell \right) }\rightarrow
A_{\left( \ell \right) }$ for $\ell \in \omega $ such that:

\begin{itemize}
\item $f_{\left( \ell +1\right) }\eta _{\left( \ell +1,\ell \right) }$ and $%
\eta _{\left( \ell +1,\ell \right) }f_{\left( \ell \right) }$ are homotopic;

\item $\mathrm{co\mathrm{lim}}_{\ell }\left( \mathrm{H}^{\bullet }\left(
f_{\left( \ell \right) }\right) \right) :\mathrm{co\mathrm{lim}}_{\ell
}{}H^{\bullet }(B_{\left( \ell \right) })\rightarrow \mathrm{co\mathrm{lim}}%
_{\ell }H^{\bullet }(A_{\left( \ell \right) })$ is an isomorphism.
\end{itemize}

Then $\mathrm{holim}\boldsymbol{A}$ and $\mathrm{ho\mathrm{lim}}\boldsymbol{B%
}$ are quasi-isomorphic.
\end{lemma}

\begin{proof}
For $\ell \in \omega $, fix a homotopy 
\begin{equation*}
f_{\left( \ell +1,\ell \right) }:f_{\left( \ell +1\right) }\eta _{\left(
\ell +1,\ell \right) }\Rightarrow \eta _{\left( \ell +1,\ell \right)
}f_{\left( \ell \right) }\text{.}
\end{equation*}

Define the map $f_{\left( \infty \right) }:\mathrm{holim}\boldsymbol{B}%
\rightarrow \mathrm{ho\mathrm{lim}}\boldsymbol{A}$ by setting%
\begin{equation*}
f_{\left( \infty \right) }^{n}\left( ze_{\ell }\right) =f_{\left( \ell
\right) }^{n}\left( z\right) e_{t_{\ell }}
\end{equation*}%
and 
\begin{equation*}
f_{\left( \infty \right) }^{n}\left( we_{\ell ,\ell +1}\right) =f_{\left(
\ell \right) }^{n+1}(w)e_{\ell ,\ell +1}+\left( -1\right) ^{n}f_{\left( \ell
+1,\ell \right) }^{n+1}\left( w\right) e_{\ell +1}\text{.}
\end{equation*}%
It is easy to verify that $f_{\left( \infty \right) }$ is indeed a map of
complexes. Furthermore, the diagram below

\begin{tikzcd}
\mathrm{H}^{\bullet }\left( \mathrm{hocolim}\boldsymbol{A}\right) \arrow[d, swap, "\mathrm{H}^{\bullet }(f_{(\infty )})"] \arrow[r] & \mathrm{colim}_{\ell }\mathrm{H}^{\bullet }(B_{\left( \ell \right)}) \arrow[d, "\mathrm{colim}_{\ell }(\mathrm{H}^{\bullet}(f_{(\ell )}))"] \\
\mathrm{H}^{\bullet }\left( \mathrm{hocolim}\boldsymbol{A}\right) \arrow[r] & \mathrm{colim}_{\ell }\mathrm{H}^{\bullet }(A_{\left( \ell \right)})
\end{tikzcd}

commutes, where the horizontal arrows are the quasi-isomorphisms from Lemma %
\ref{Lemma:cohomology-colimit}. Hence $f_{(\infty )}$ is a quasi-isomorphism
of complexes over $\mathcal{B}$.
\end{proof}

\begin{lemma}
\label{Lemma:qi-2}Suppose that $\boldsymbol{A}$ and $\boldsymbol{B}$ are
inductive sequences of complexes over $\mathcal{B}$. Suppose that there
exists an increasing sequence $\left( t_{\ell }\right) $ in $\omega $ and
maps $f_{\left( \ell \right) }:B_{\left( \ell \right) }\rightarrow A_{\left(
t_{\ell }\right) }$ for $\ell \in \omega $ such that:

\begin{itemize}
\item $f_{\left( \ell +1\right) }\eta _{\left( \ell +1,\ell \right) }$ and $%
\eta _{\left( t_{\ell +1},t_{\ell }\right) }f_{\left( \ell \right) }$ are
cohomologous;

\item $\mathrm{co\mathrm{lim}}_{\ell }\left( \mathrm{H}^{\bullet }\left(
f_{\left( \ell \right) }\right) \right) :\mathrm{co\mathrm{lim}}_{\ell
}H^{\bullet }(B_{\left( \ell \right) })\rightarrow \mathrm{co\mathrm{lim}}%
_{\ell }H^{\bullet }(A_{\left( \ell \right) })$ is an isomorphism.
\end{itemize}

Then $\mathrm{holim}\boldsymbol{A}$ and $\mathrm{ho\mathrm{lim}}\boldsymbol{B%
}$ are quasi-isomorphic.
\end{lemma}

\begin{proof}
Without loss of generality, we can assume that $t_{\ell }=\ell $ for $\ell
\in \omega $. Using the fact that quasi-isomorphisms in $\mathbf{K}\left( 
\mathcal{B}\right) $ are a multiplicative system, one can easily define an
inductive sequence $\boldsymbol{C}$ of complexes over $\mathcal{B}$, and a
sequence $\left( g_{\left( \ell \right) }\right) $ of maps $g_{\ell
}:A_{\left( \ell \right) }\rightarrow C_{\left( \ell \right) }$ such that $%
\left( g_{\left( \ell \right) }\right) $ and $\left( g_{\left( \ell \right)
}f_{\left( \ell \right) }\right) $ satisfy the assumptions of Lemma \ref%
{Lemma:qi-1}. Thus we obtain that $\mathrm{hoco\mathrm{lim}}\boldsymbol{A}$
and \textrm{hocolim}$\boldsymbol{C}$, and, respectively, $\mathrm{hoco%
\mathrm{lim}}\boldsymbol{B}$ and $\mathrm{hoco\mathrm{lim}}\boldsymbol{C}$
are quasi-isomorphic by applying Lemma \ref{Lemma:qi-1} twice.
\end{proof}

Suppose now that $\mathcal{A}$ is a countably complete abelian category. Let 
\begin{equation*}
\boldsymbol{A}=\left( A^{\left( n\right) },p^{\left( n,m\right) }:A^{\left(
m\right) }\rightarrow A^{\left( n\right) }\right)
\end{equation*}%
be a tower of complexes over $\mathcal{A}$. As in \cite[Section 17]%
{mardesic_strong_2000} or \cite%
{cordier_limites_1987,cordier_homotopy_1997,cordier_homologie_1987}, one
defines the \emph{homotopy limit }$\mathrm{ho\mathrm{\mathrm{lim}}}%
\boldsymbol{A}$ of $\boldsymbol{A}$ to be the complex over $\mathcal{A}$
defined by setting%
\begin{equation*}
\left( \mathrm{ho\mathrm{\mathrm{lim}}}\boldsymbol{A}\right) _{n}=\mathrm{%
\mathrm{Ker}}\left( \eta :\prod_{m\in \omega }A_{n}^{\left( m\right) }\oplus
\prod_{m_{0}\leq m_{1}}A_{n+1}^{\left( m_{0}\right) }\rightarrow
\prod_{m_{0}\leq m_{1}\leq m_{2}}A_{n+1}^{\left( m_{0}\right) }\right) \text{%
.}
\end{equation*}%
Here, $\eta $ is defined by%
\begin{equation*}
\pi _{m_{0},m_{1},m_{2}}\left( x\right) =x_{m_{0},m_{1}}+p^{\left(
m_{0},m_{1}\right) }(x_{m_{1},m_{2}})-x_{m_{0},m_{2}}
\end{equation*}%
for each element $x$. The differential 
\begin{equation*}
\partial _{n}:\left( \mathrm{ho\mathrm{\mathrm{lim}}}\boldsymbol{A}\right)
_{n+1}\rightarrow \left( \mathrm{ho\mathrm{\mathrm{lim}}}\boldsymbol{A}%
\right) _{n}
\end{equation*}%
is defined by setting%
\begin{equation*}
\left( \partial _{n}x\right) _{m}=\partial _{n}\left( x_{m}\right)
\end{equation*}%
an%
\begin{equation*}
\left( \partial _{n}x\right) _{m_{0},m_{1}}=\partial _{n+1}\left(
x_{m_{0},m_{1}}\right) +\left( -1\right) ^{n}(p_{n}^{\left(
m_{0},m_{1}\right) }(x_{m_{1}})-x_{m_{0}})\text{.}
\end{equation*}

It is manifest that the notions of homotopy limit and homotopy colimit are
dual to each other. Let $\mathcal{M}$ be a countably complete abelian
category, and $X$ be an object of $\mathcal{M}$. Suppose that $\boldsymbol{A}%
=\left( A_{\left( m\right) }\right) $ is an inductive sequence of complexes
of countable abelian groups. Adopting the notation of Section \ref%
{Subsection:derive-hom}, we obtain a tower 
\begin{equation*}
\mathrm{Hom}_{\mathcal{M}}\left( \boldsymbol{A},X\right) :=\left( \mathrm{Hom%
}\left( A_{\left( n\right) },X\right) \right)
\end{equation*}%
of complexes in $\mathcal{M}$. Furthermore, since $\mathrm{Hom}_{\mathcal{M}%
}\left( -,X\right) $ is continuous, we have that%
\begin{equation*}
\mathrm{ho\mathrm{lim}}\left( \mathrm{Hom}\left( \boldsymbol{A},X\right)
\right) \cong \mathrm{Hom}\left( \mathrm{hoco\mathrm{lim}}\boldsymbol{A}%
,X\right) \text{.}
\end{equation*}%
Furthermore, by Corollary \ref{Corollary:UCT} we have that if $\mathrm{hoco%
\mathrm{lim}}\boldsymbol{A}$ and $\mathrm{hoco\mathrm{lim}}\boldsymbol{B}$
are quasi-isomorphic, then $\mathrm{ho\mathrm{lim}}\left( \mathrm{Hom}\left( 
\boldsymbol{A},X\right) \right) $ and $\mathrm{ho\mathrm{lim}}\left( \mathrm{%
Hom}\left( \boldsymbol{B},X\right) \right) $ are quasi-isomorphic.

\section{Universal Coefficient Theorem for homology and cohomology\label%
{Section:topology}}

In this section we present a definition of Steenrod homology \cite%
{steenrod_regular_1940,steenrod_regular_1941} for compact Polish spaces and 
\v{C}ech cohomology \cite{eilenberg_foundations_1952} for (homotopy)
polyhedra with coefficients in an arbitrary countably complete abelian
category. The general Universal Coefficient Theorem is then applied to
obtain Universal Coefficient Sequences expressing Steenrod homology in terms
of \v{C}ech cohomology and vice versa.

\subsection{Simplicial complexes\label{Subsection:simplicial-complexes}}

We begin with recalling some fundamental notions concerning simplicial
complexes and their homology and cohomology as can be found, for example, in 
\cite{eilenberg_foundations_1952}. A \emph{simplicial complex} $K$ is a
family of nonempty finite sets that is closed downwards, i.e., $\sigma
\subseteq \tau \in K\Rightarrow \sigma \in K$. A \emph{simplex }of\emph{\ }$%
K $ is any element $\sigma \in K$. If $\sigma ,\sigma ^{\prime }$ are
simplices of $K$ such that $\sigma \subseteq \sigma ^{\prime }$, then $%
\sigma $ is a face of $\sigma ^{\prime }$. A \emph{vertex of $K$} is any
element $v$ of $\mathrm{dom}(K):=\bigcup K$. Let $K,L$ be two simplicial
complexes. A \emph{simplicial map} $f\colon K\rightarrow L$ is any function $%
f\colon \mathrm{dom}(K)\rightarrow \mathrm{dom}(L)$ so that $%
\{f(v_{0}),\ldots ,f(v_{n})\}\ $a simplex of $L$ for every simplex $%
\{v_{0},\ldots ,v_{n}\}\ $of $K$. Two simplicial maps $f,f^{\prime
}:K\rightarrow L$ are \emph{contiguous }if $\left\{ f\left( v\right)
,f^{\prime }\left( v\right) \right\} $ is a simplex for every $v\in K$. A
simplicial complex is \emph{finite} if it has finitely many vertices, \emph{%
countable }if it has countably many vertices, and \emph{locally finite }if
every vertex belongs to finitely many simplices.

One can associate with a countable locally finite simplicial complex $K$ a
locally compact Polish space, called its \emph{geometric realization}, as
follows. Let $\Xi $ be the separable Hilbert space with orthonormal basis $%
\left( e_{v}\right) _{v\in \mathrm{\mathrm{dom}}(K)}$ indexed by the set of
vertices of $K$. For each simplex $\sigma =\left\{ v_{0},\ldots
,v_{n}\right\} $ of $K$ define%
\begin{equation*}
\left\vert \sigma \right\vert =\left\{ t_{0}e_{v_{0}}+\cdots
+t_{n}e_{v_{n}}:t_{0},\ldots ,t_{n}\in \left[ 0,1\right] ,t_{0}+\cdots
+t_{n}=1\right\} \subseteq \Xi \text{.}
\end{equation*}%
Then $\left\vert K\right\vert $ is defined to be the closed subspace of $\Xi 
$ obtained as the union of $\left\vert \sigma \right\vert $ where $\sigma $
ranges among the simplices of $K$. A \emph{polyhedron }is a locally compact
Polish space that is obtained in this fashion from some simplicial complex 
\cite[Chapter 1]{munkres_elements_1984}; see also \cite[Section II.3,
Proposition 3.6]{lundell_topology_1969}.

We now recall the definition of the classical homology invariants of a
simplicial complex. We adopt the notation and terminology from \cite[Chapter
VI]{eilenberg_foundations_1952}. For $n\geq 0$, an \emph{elementary} $n$-%
\emph{chain} is any tuple $(v_{0},\ldots ,v_{n})\in \mathrm{dom}(K)^{n+1}$
with $\{v_{0},\ldots ,v_{n}\}$ a simplex of $K$. The \emph{chain complex} $%
C_{\bullet }\left( K\right) $ of $K$ is defined by letting $C_{n}(K)$ be the
free abelian group generated by the set of elementary $n$-chains for $n\geq
0 $ \cite[Definition 2.3]{eilenberg_foundations_1952}. Elements of $C_{n}(K)$
are called (\emph{ordered}) $n$-\emph{chains\ }of $K$. The differential $%
\partial :C_{n}(K)\rightarrow C_{n-1}(K)$ is defined by 
\begin{equation*}
\partial (v_{0},\ldots ,v_{n})=\sum_{i=0}^{n}\left( -1\right)
^{i}(v_{0},\ldots ,\hat{v}_{i},\ldots ,v_{n})\text{,}
\end{equation*}%
where $\hat{v}_{i}$ denotes the omission of $v_{i}$. The cochain complex $%
C^{\bullet }\left( K\right) $ is defined by setting $C^{n}\left( K\right) =%
\mathrm{Hom}\left( C_{n}\left( K\right) ,\mathbb{Z}\right) $ for $n\geq 0$.
Thus, we have that $C^{n}\left( K\right) $ is the product of copies of $%
\mathbb{Z}$ indexed by elementary $n$-chains for $n\geq 0$. The
codifferential $\delta :C^{n}\left( K\right) \rightarrow C^{n+1}\left(
K\right) $ is defined by setting%
\begin{equation*}
\left( \delta f\right) \left( v_{0},\ldots ,v_{n+1}\right)
=\sum_{i=0}^{n}\left( -1\right) ^{i}f(v_{0},\ldots ,\hat{v}_{i},\ldots
,v_{n+1})
\end{equation*}%
for $f\in C^{n}\left( K\right) $ and an elementary $\left( n+1\right) $%
-chain $\left( v_{0},\ldots ,v_{n+1}\right) $. A simplicial map $%
f:K\rightarrow L$ induces a morphism of complexes $C_{\bullet }\left(
K\right) \rightarrow C_{\bullet }\left( L\right) $, in such a way that
contiguous simplicial maps induce homotopic morphisms.

The \emph{barycentric subdivision }$\mathrm{Sd}\left( K\right) $ of $K$ is
the simplicial complex with $\mathrm{dom}\left( \mathrm{Sd}\left( K\right)
\right) $ equal to the set of nonempty simplices of $K$. A simplex in $%
\mathrm{Sd}\left( K\right) $ is a set $\left\{ \sigma _{0},\ldots ,\sigma
_{n}\right\} $ of simplices of $K$ that is linearly ordered by inclusion. A
vertex selection map for $K$ is a function $f:\mathrm{\mathrm{dom}}\left( 
\mathrm{Sd}\left( K\right) \right) \rightarrow \mathrm{\mathrm{dom}}\left(
K\right) $ such that $f\left( \sigma \right) \in \sigma $ for every $\sigma
\in \mathrm{\mathrm{dom}}\left( \mathrm{Sd}\left( K\right) \right) $. This
can be regarded as a simplicial map $\mathrm{Sd}\left( K\right) \rightarrow
K $. It induces a quasi-isomorphism $C_{\bullet }\left( \mathrm{Sd}\left(
K\right) \right) \rightarrow C_{\bullet }\left( K\right) $, which does not
depend on the choice of the vertex selection map (as any two such maps are
contiguous).

\subsection{Homology of compact Polish spaces\label{Section:homology-spaces}}

Let $X$ be a compact Polish space. A \emph{cover} of $X$ is a family $%
\mathcal{U}=\left( U_{i}^{\mathcal{U}}\right) _{i\in \omega }$ of open
subsets of $X$ such that $X$ is the union of $\left\{ U_{i}^{\mathcal{U}%
}:i\in \omega \right\} $. The \emph{\v{C}ech nerve} $N(\mathcal{U})$ of a
cover $\mathcal{U}$ of $X$ is the simplicial complex with 
\begin{equation*}
\mathrm{dom}(N(\mathcal{U}))=\mathrm{supp}(\mathcal{U}):=\left\{ i\in \omega
:U_{i}^{\mathcal{U}}\neq \varnothing \right\}
\end{equation*}%
and $\sigma =\left\{ v_{0},\ldots ,v_{n}\right\} $ a simplex of $N\left( 
\mathcal{U}\right) $ if and only if 
\begin{equation*}
U_{\sigma }^{\mathcal{U}}:=U_{v_{0}}^{\mathcal{U}}\cap \cdots \cap
U_{v_{n}}^{\mathcal{U}}\neq \varnothing \text{.}
\end{equation*}%
We say that $\mathcal{U}$ is finite if $N\left( \mathcal{U}\right) $ is
finite. The \emph{Vietoris nerve }$V\left( \mathcal{U}\right) $ is the
simplicial complex with 
\begin{equation*}
\mathrm{\mathrm{dom}}\left( V\left( \mathcal{U}\right) \right) =X
\end{equation*}%
and $\sigma =\left\{ x_{0},\ldots ,x_{n}\right\} $ is a simplex if and only
if there exists $i\in \omega $ such that $\sigma \subseteq U_{i}^{\mathcal{U}%
}$. A \emph{point selection map} for $\mathcal{U}$ is a function $g:\omega
\rightarrow X$ such that $g\left( i\right) \in U_{i}^{\mathcal{U}}$ for
every $i\in X$. This defines a simplicial map \textrm{Sd}$\left( N\left( 
\mathcal{U}\right) \right) \rightarrow V\left( \mathcal{U}\right) $,
inducing a quasi-isomorphism $C^{\bullet }\left( V\left( \mathcal{U}\right)
\right) \rightarrow C^{\bullet }\left( \mathrm{Sd}\left( N\left( \mathcal{U}%
\right) \right) \right) $ which does not depend on the choice of the point
selection map (as any two such maps are contiguous).

A cover $\mathcal{U}$ of $X$ is a \emph{refinement} of a cover $\mathcal{V}$
of $X$ if for every $i\in \mathrm{\mathrm{\mathrm{supp}}}(\mathcal{U})$
there exists $f\left( i\right) \in \mathrm{\mathrm{\mathrm{supp}}}(\mathcal{V%
})$ such that $U_{i}^{\mathcal{U}}\subseteq U_{j}^{\mathcal{V}}$. In this
case we set $\mathcal{V}\leq \mathcal{U}$ and call $f$ a \emph{refinement map%
}. We have that $f$ is a simplicial map $N\left( \mathcal{U}\right)
\rightarrow N\left( \mathcal{V}\right) $, while the identity is a simplicial
map $V\left( \mathcal{U}\right) \rightarrow V\left( \mathcal{V}\right) $.

Let $X$ be a compact Polish space. The relation $\,\leq \,$ renders the set $%
\mathtt{cov}(X)$ of \emph{finite} covers of $X$ an upward directed ordering
of countable cofinality. A \emph{cofinal sequence }of finite open covers for 
$X$ is a cofinal increasing sequence $\mathcal{U}=\left( \mathcal{U}%
_{m}\right) _{m\in \omega }$ in $\mathtt{cov}(X)$. This together with a
choice of refinement maps $\mathcal{U}_{m+1}\rightarrow \mathcal{U}_{m}$
yield an inductive sequence of cochain complexes%
\begin{equation*}
\underline{\boldsymbol{C}}^{\bullet }\left( X\right) :=\left( C^{\bullet
}\left( V\left( \mathcal{U}_{m}\right) \right) \right) _{m\in \omega }
\end{equation*}%
and an inductive sequence of cochain complexes of free finitely-generated
abelian groups%
\begin{equation*}
\boldsymbol{C}^{\bullet }\left( X\right) :=\left( C^{\bullet }\left( N\left( 
\mathcal{U}_{m}\right) \right) \right) _{m\in \omega }\text{.}
\end{equation*}%
We then define%
\begin{equation*}
\underline{C}^{\bullet }\left( X\right) :=\mathrm{hoco\mathrm{lim}}%
\underline{\boldsymbol{C}}^{\bullet }\left( X\right)
\end{equation*}%
and%
\begin{equation*}
C^{\bullet }\left( X\right) :=\mathrm{hoco\mathrm{lim}}\boldsymbol{C}%
^{\bullet }\left( X\right)
\end{equation*}%
The assignment $X\mapsto \underline{\boldsymbol{C}}^{\bullet }\left(
X\right) $ is a contraviariant functor from the category $\mathbf{PK}$ of
compact Polish spaces to $\mathrm{Pro}_{\omega }\left( \mathbf{K}^{-}\left( 
\mathbf{Ab}\right) \right) $. Thus, we have that $X\mapsto \underline{C}%
^{\bullet }\left( X\right) :=\mathrm{hoco\mathrm{lim}}\underline{\boldsymbol{%
C}}^{\bullet }\left( X\right) $ can be seen as a functor $\mathbf{PK}^{%
\mathrm{op}}\rightarrow \mathbf{D}^{-}\left( \mathbf{Ab}\right) $.

The choice of cover selection maps for $\mathcal{U}_{m}$ and of vertex
selection maps for $N\left( \mathcal{U}_{m}\right) $ yield quasi-isomorphisms%
\begin{equation*}
C^{\bullet }\left( V\left( \mathcal{U}_{m}\right) \right) \rightarrow
C^{\bullet }\left( N\left( \mathcal{U}_{m}\right) \right)
\end{equation*}%
for $m\in \omega $ that satisfy the assumptions of Lemma \ref{Lemma:qi-2};
see \cite[Section 5]{dowker_homology_1952}. Hence, we have that $C^{\bullet
}\left( X\right) $ and $\underline{C}^{\bullet }\left( X\right) $ are
quasi-isomorphic. Thus, $X\mapsto C^{\bullet }\left( X\right) $ yields a
functor from $\mathbf{PK}$ to $\mathbf{D}^{-}\left( \mathbf{Ab}_{\aleph
_{0}}\right) $, where $\mathbf{Ab}_{\aleph _{0}}$ is the category of
countable abelian groups.

Let now $\mathcal{M}$ be a countably complete abelian category, and let $G$
be an object of $\mathcal{M}$. Then by Corollary \ref{Corollary:UCT} we have
that setting 
\begin{equation*}
C_{\bullet }\left( X;G\right) :=\mathrm{Hom}\left( C^{\bullet }\left(
X\right) ,G\right) \cong \mathrm{ho\mathrm{lim}}\left( \mathrm{Hom}\left(
C_{\bullet }\left( \mathcal{U}_{m}\right) ,G\right) \right)
\end{equation*}%
yields a functor $\mathbf{PK}\rightarrow \mathbf{D}^{+}\left( \mathcal{M}%
\right) $. We then define the \emph{Steenrod homology} $\mathrm{H}_{\bullet
}\left( X;G\right) $ of $X$ with coefficients in $G$ to be the homology of
the chain complex $C_{\bullet }\left( X;G\right) $ over $\mathcal{M}$. We
also define $\mathrm{H}^{\bullet }\left( X\right) $ to be the cohomology of
the cochain complex of countable abelian groups $C^{\bullet }\left( X\right) 
$. As a particular instance of Theorem \ref{Theorem:UCT} we have the
following.

\begin{theorem}
\label{Theorem:UCT-homology-spaces}Suppose that $X$ is a compact Polish
space. Let $\mathcal{M}$ be a countably complete abelian category, and $G$
be an object of $\mathcal{M}$. Then we have a natural short exact sequence
of graded objects of $\mathcal{M}$%
\begin{equation*}
0\rightarrow \mathrm{Ext}\left( \mathrm{H}^{\bullet +1}\left( X\right)
,G\right) \rightarrow \mathrm{H}_{\bullet }\left( X;G\right) \rightarrow 
\mathrm{Hom}\left( \mathrm{H}^{\bullet }\left( X\right) ,G\right)
\rightarrow 0
\end{equation*}%
which splits (unnaturally).
\end{theorem}

Theorem \ref{Theorem:UCT-homology-spaces} in particular entails the
homotopy-invariance of $\mathrm{H}_{\bullet }\left( -;G\right) $, in view of
the homotopy-invariance of integral cohomology. More generally, one can
consider definable Steenrod homology of \emph{pairs of compact Polish spaces 
}and obtain a corresponding Universal Coefficient Sequence.

\subsection{Cohomology of (homotopy) polyhedra\label%
{Subsection:cohomology-polyhedra}}

In this section we assume all the simplicial complexes to be \emph{countable 
}and \emph{locally finite}. Recall that a polyhedron is a topological space $%
X$ that is the geometric realization $\left\vert K\right\vert $ of a
simplicial complex $K$. In particular, a polyhedron is a locally compact
Polish space. If $\left\vert K\right\vert $ is a polyhedron, then one
defines its corresponding chain complex $C_{\bullet }^{\mathrm{Pol}%
}(\left\vert K\right\vert )$ to be the free countable chain complex $%
C_{\bullet }(K)$ associated with the simplicial complex $K$. This defines a
functor $X\mapsto C_{\bullet }^{\mathrm{Pol}}\left( X\right) $ from the
category $\mathbf{Pol}$\textbf{\ }of polyhedra to $\mathbf{K}^{+}\left( 
\mathbf{Ab}_{\aleph _{0}}\right) $; see \cite[Sections 15--18]%
{munkres_elements_1984}. One then defines $\mathrm{H}_{\bullet }\left(
X\right) $ to be the corresponding homology, which is naturally isomorphic
to the homology defined in the previous section in the particular case when $%
X$ is compact.

Let $\mathcal{M}$ be a countably complete abelian category, and let $G$ be
an object of $\mathcal{M}$. We define, for a polyhedron $Y$,%
\begin{equation*}
C_{\mathrm{Pol}}^{\bullet }\left( Y;G\right) :=\mathrm{Hom}\left( C_{\bullet
}^{\mathrm{Pol}}\left( Y\right) ,G\right) \text{.}
\end{equation*}%
This gives a functor from $\mathbf{Pol}^{\mathrm{op}}$ to $\mathbf{D}%
^{-}\left( \mathcal{M}\right) $. We define $\mathrm{H}^{\bullet }\left(
Y;G\right) $ to be the corresponding cohomology. As a particular instance of
Theorem \ref{Theorem:UCT}, we have the following.

\begin{theorem}
\label{Theorem:UCT-cohomology-polyhedra}Suppose that $Y$ is a countable,
locally finite polyhedron. Let $\mathcal{M}$ be a countable complete abelian
category, and $G$ be an object of $\mathcal{M}$. Then we have a natural
short exact sequence of graded objects of $\mathcal{M}$%
\begin{equation*}
0\rightarrow \mathrm{Ext}\left( \mathrm{H}_{\bullet -1}\left( Y\right)
,G\right) \rightarrow \mathrm{H}^{\bullet }\left( Y;G\right) \rightarrow 
\mathrm{Hom}\left( \mathrm{H}_{\bullet }\left( Y\right) ,G\right)
\rightarrow 0
\end{equation*}%
which splits (unnaturally).
\end{theorem}

More generally, one can consider definable \v{C}ech cohomology of \emph{%
pairs of polyhedra }and obtain a corresponding Universal Coefficient
Sequence.

Suppose that $Y=\left\vert K\right\vert $ is a polyhedron. A \emph{%
polyhedral cofiltration }of $Y$ is an increasing sequence $\left( \left\vert
K_{m}\right\vert \right) _{m\in \omega }$ of compact subspaces of $Y$ such
that $\left( K_{m}\right) _{m\in \omega }$ is an increasing subsequence of 
\emph{finite} subcomplexes of $K$ with union equal to $K$. One can then
consider the corresponding\emph{\ }inductive sequence of chain complexes of
free finitely-generated abelian groups $\left( C_{\bullet }\left( \left\vert
K_{m}\right\vert \right) \right) _{m\in \omega }$. Then we have that $%
C_{\bullet }\left( \left\vert K\right\vert \right) $ is quasi-isomorphic to $%
\mathrm{colim}_{m}C_{\bullet }(\left\vert K_{m}\right\vert )$, which in turn
is quasi-isomorphic to $\mathrm{hoco\mathrm{lim}}_{m}C_{\bullet }\left(
\left\vert K_{m}\right\vert \right) $. Thus, we have that $C^{\bullet
}\left( Y;G\right) $ is quasi-isomorphic to $\mathrm{ho\mathrm{lim}}%
_{m}C^{\bullet }\left( \left\vert K_{m}\right\vert ;G\right) $.

By homotopy invariance, one can extend the definition of homology and
cohomology to the class of \emph{homotopy polyhedra}, namely spaces that are 
\emph{homotopy equivalent }to a polyhedron. Indeed, if $X$ is homotopy
equivalent to a polyhedron $\left\vert K_{X}\right\vert $, as witnessed by
continuous maps $s_{X}:X\rightarrow \left\vert K_{X}\right\vert $ and $%
t_{X}:\left\vert K_{X}\right\vert \rightarrow X$, then one can set $\mathrm{H%
}^{\bullet }\left( X;G\right) :=\mathrm{H}^{\bullet }\left( \left\vert
K_{X}\right\vert ;G\right) $ and $\mathrm{H}_{\bullet }\left( X\right) :=%
\mathrm{H}_{\bullet }\left( \left\vert K_{X}\right\vert \right) $. If $%
\alpha :X\rightarrow Y$ is a continuous map, then the corresponding
morphisms $\mathrm{H}^{\bullet }\left( Y;G\right) \rightarrow \mathrm{H}%
^{\bullet }\left( X;G\right) $ and $\mathrm{H}_{\bullet }\left( X\right)
\rightarrow \mathrm{H}_{\bullet }\left( Y\right) $ are the ones associated
with the continuous map $s_{Y}\circ \alpha \circ t_{X}:\left\vert
K_{X}\right\vert \rightarrow \left\vert K_{Y}\right\vert $.

The class of homotopy polyhedra contains all \emph{countable CW complexes} 
\cite[Section\textrm{\ }1.5]{arkowitz_introduction_2011}. A CW complex is
countable if it is obtained by attaching countably many cells. Every
polyhedron is a countable simplicial complex. Conversely, every countable CW
complex is a homotopy polyhedron by \cite[Theorem 1]{milnor_spaces_1959}, 
\cite[Section IV.6, Theorem 6.1]{lundell_topology_1969}. 
Similar considerations apply to absolute neighborhood retract (ANR) \cite[%
Section I.3]{mardesic_shape_1982}: every polyhedron is an ANR and,
conversely, every ANR is a homotopy polyhedron \cite[Section I.4, Theorem 5]%
{mardesic_shape_1982}.

We also have that a \emph{second countable paracompact space with a good
cover} is a homotopy polyhedron. Recall that a second countable topological
spaces is \emph{paracompact} if every cover has a countable, locally finite
refinement. A \emph{good cover} of a paracompact space $X$ is a countable,
locally finite cover $\mathcal{U}$ of $X$ such that, for every $\sigma \in N(%
\mathcal{U})$, $U_{\sigma }^{\mathcal{U}}$ is contractible. In this case,
one has that $X$ is homotopy equivalent to the polyhedron $\left\vert N(%
\mathcal{U})\right\vert $ \cite[Corollary 4G.3]{hatcher_algebraic_2002}.

The class of second countable paracompact spaces with a good cover includes
all locally compact Polish spaces admitting a basis that (1) is closed under
intersections, and (2) consists of precompact contractible open sets. In
particular, all second countable Riemannian manifolds have this property;
see \cite[Chapter 11]{spivak_comprehensive_1970}, \cite[Theorem 5.1]%
{bott_differential_1982}, \cite[proof of Theorem 89]%
{petersen_riemannian_2006}, \cite[Section 3.4]{do_carmo_riemannian_1992}, 
\cite[Remark after Lemma 10.3]{milnor_morse_1963}.

\section{Complexity theory and the left heart of Polish abelian groups\label%
{Section:complexity}}

In this section we recall some notions from topology and descriptive set
theory, including Polish groups, groups with a Polish cover \cite%
{bergfalk_definable_2024,bergfalk_definable_2024-1}, and complexity of
equivalence relations. In particular, we explain how groups with a Polish
cover provide a concrete description of the left heart of Polish groups, as
shown in \cite{lupini_looking_2024}. Standard sources for descriptive set
theory and the theory of Polish groups and their actions include \cite%
{kechris_classical_1995,gao_invariant_2009,becker_descriptive_1996}.

\subsection{Ordinal notation\label{Subsection:notation}}

Notationally, we let $\omega $ be the first infinite ordinal, and $\mathbb{N}
$ be the set of strictly positive integers. For each countable limit ordinal 
$\lambda $ we fix a strictly increasing sequence $\left( \lambda _{n}\right)
_{n\in \omega }$ of \emph{successor }ordinals with $\lambda =\mathrm{sup}%
_{n}{}\lambda _{n}$. For a zero or successor ordinal $\alpha $, define $%
\alpha _{n}=\alpha $ for $n\in \omega $. We also set $\alpha _{n}=-1$ when $%
\alpha =-1$ by convention.

We define $\omega _{1}[1/2]$ to be the ordered set containing the set $%
\omega _{1}$ of all countable ordinals together with $\alpha +1/2$ for $%
\alpha \in \omega _{1}$ subject to the relations%
\begin{equation*}
\alpha <\alpha +\frac{1}{2}<\alpha +1\text{.}
\end{equation*}%
For a successor ordinal $\alpha $, we define $\alpha -1$ to be its immediate
predecessor and%
\begin{equation*}
\alpha -\frac{1}{2}=\left( \alpha -1\right) +\frac{1}{2}\text{.}
\end{equation*}%
We also set $\alpha -1=-1$ for $\alpha =0$. Let now $\omega _{1}^{\mathbf{Pol%
}}$ be the ordered set containing $\omega _{1}[1/2]$ together with $\lambda +%
\frac{1}{2}+\varepsilon $ for a limit or zero ordinal $\lambda $, subject to
the relations%
\begin{equation*}
\lambda +\frac{1}{2}<\lambda +\varepsilon <\lambda +1\text{.}
\end{equation*}

\subsection{Invariant complexity theory}

Descriptive set theory \cite{kechris_classical_1995} studies the properties
of \textquotedblleft definable\textquotedblright\ subset of
\textquotedblleft well-behaved\textquotedblright\ topological spaces, such
as Polish spaces (second countable spaces whose topology is induced by a
complete metric). Within descriptive set theory, invariant complexity theory 
\cite{gao_invariant_2009} focuses on the study of the complexity of \emph{%
classification problems }in mathematics. In this context, a classification
problem is regarded as an equivalence relation $E$ on a Polish space $X$,
which is a setting that includes virtually all classification problems
arising in practice in mathematics. Two such classification problems $\left(
X,E\right) $ and $\left( Y,F\right) $ are compared via the notion of Borel
reduction, which is an injection $X/E\rightarrow Y/F$ between the quotient
spaces that is \emph{Borel-definable}, in the sense that it is induced by a
Borel function $X\rightarrow Y$. In this case, we say that $E$ is\emph{\
Borel reducible }to $F$. This captures the idea that the classification
problem represented by $\left( X,E\right) $ is at most as complicated as the
one represented by $\left( Y,F\right) $. If $E$ is Borel reducible to $F$
and vice versa, then $E$ and $F$ are called \emph{Borel bireducible}.

Distinguished classification problems serve as \emph{benchmarks }to measure
the complexity of all the other ones. Such benchmarks include the relation $%
=_{\mathbb{R}}$ of equality on $\mathbb{R}$ (or, equivalently, any
uncountable Polish space), and the relation $E_{0}$ of \emph{tail
equivalence }on the space $\mathcal{C}:=\left\{ 0,1\right\} ^{\mathbb{N}}$
of \emph{binary sequences}.

By a complexity class, we mean a function $X\mapsto \Gamma \left( X\right) $
that assigns to each Polish space $X$ a collection $\Gamma \left( X\right) $
of Borel subsets of $X$, in such a way that $f^{-1}\left( A\right) \in
\Gamma \left( X\right) $ for every continuous function $f:X\rightarrow Y$
between Polish spaces and $A\in \Gamma \left( Y\right) $. We define its dual
class $\check{\Gamma}\left( X\right) $ by setting $A\in \check{\Gamma}\left(
X\right) \Leftrightarrow X\setminus A\in \Gamma \left( X\right) $. We say
that $\Gamma $ is self-dual if $\check{\Gamma}=\Gamma $. For $\alpha <\omega
_{1}$, besides $\boldsymbol{\Pi }_{\alpha }^{0}$ and $\boldsymbol{\Sigma }%
_{\alpha }^{0}$, one can also consider the complexity class $D(\boldsymbol{%
\Pi }_{\alpha }^{0})$ comprising sets that are the set-theoretic difference
of sets in $\boldsymbol{\Pi }_{\alpha }^{0}$. If $\Gamma $ is a complexity
class that is not self-dual, and $A$ is a subset of a Polish space $X$, then
we say that $\Gamma $ is the complexity class of $A$ if $A\in \Gamma \left(
X\right) $ and $X\setminus A\notin \Gamma \left( X\right) $.

If $E$ is an equivalence relation on a Polish space $X$, then we say that $E$
is \emph{potentially} $\Gamma $ if it is Borel reducible to an equivalence
relation $F$ on a Polish space $Y$ such that $F\in \Gamma \left( Y\times
Y\right) $. When $\Gamma $ is not self-dual, we say that $\Gamma $ is \emph{%
the} potential Borel complexity class of $E$ if $E$ is potentially $\Gamma $
and $E$ is not potentially $\check{\Gamma}$. If $X/E$ is $\boldsymbol{\Sigma 
}_{1}^{1}$-definable set, then we say that $E$ is $\Gamma $-definable if and
only if $E$ is potentially $\Gamma $, and $\Gamma $ is the complexity class
of $X/E$ if and only if $\Gamma $ is the potential complexity class of $E$.

We define complexity classes $\Gamma _{\alpha }$ for every $\alpha \in
\omega _{1}^{\mathbf{Pol}}$ as follows. For $\lambda <\omega _{1}$ zero or
limit, $n<\omega $, and $i\in \left\{ 0,1/2,1/2+\varepsilon \right\} $:

\begin{equation*}
\Gamma _{\lambda +n+i}:=\left\{ 
\begin{array}{ll}
\boldsymbol{\Pi }_{1+\lambda }^{0} & \text{for }n=0\text{ and }i=0\text{;}
\\ 
\boldsymbol{\Sigma }_{1+\lambda +1}^{0} & \text{for }n=0\text{ and }i=\frac{1%
}{2}\text{;} \\ 
D(\boldsymbol{\Pi }_{1+\lambda +1}^{0}) & \text{for }n=0\text{ and }i=\frac{1%
}{2}+\varepsilon \text{;} \\ 
D(\boldsymbol{\Pi }_{1+\lambda +n+1}^{0}) & \text{for }n\geq 1\text{ and }i=%
\frac{1}{2}\text{.} \\ 
\boldsymbol{\Pi }_{1+\lambda +n+1}^{0} & \text{for }n\geq 1\text{ and }i=0%
\text{;}%
\end{array}%
\right. 
\end{equation*}%
It is easily verified that $\left( \Gamma _{\alpha }\right) $ is an
increasing family of complexity classes.

\subsection{Borel-definable sets}

Recall that a $\boldsymbol{\Sigma }_{1}^{1}$-definable\emph{\ set }is a set $%
X$ explicitly presented as the quotient of a Polish space $\hat{X}$ by an
analytic (i.e., $\mathbf{\Sigma }_{1}^{1}$) equivalence relation $E$ \cite[%
Section 14]{kechris_classical_1995}. When $E$ is \emph{Borel} and\emph{\
idealistic }\cite[Section 3]{bergfalk_definable_2024-1}, then we say that $X$
is a\emph{\ Borel-definable set}. For example, if $G$ is a Polish group
acting continuously on a Polish space $\hat{X}$, then the corresponding
orbit equivalence relation $E_{G}^{\hat{X}}$ is analytic and idealistic, and 
$\hat{X}/E_{G}^{\hat{X}}$ is a $\mathbf{\Sigma }_{1}^{1}$-definable set. If
furthermore $E_{G}^{\hat{X}}$ is Borel, then $\hat{X}/E_{G}^{\hat{X}}$ is a
Borel-definable set. A function $f:\hat{X}/E\rightarrow \hat{Y}/F$ between $%
\mathbf{\Sigma }_{1}^{1}$-definable sets is Borel-definable if it is induced
by a Borel function $\varphi :\hat{X}\rightarrow \hat{Y}$, which means that $%
f([x]_{E})=[\varphi \left( x\right) ]_{F}$ for every $x\in \hat{X}$.
Borel-definable functions are the morphisms in the category of $\boldsymbol{%
\Sigma }_{1}^{1}$-definable sets. It is proved in \cite[Lemma 3.7]%
{kechris_borel_2016} that if $X=\hat{X}/E$ is a $\boldsymbol{\Sigma }_{1}^{1}
$-definable set such that $E$ is idealistic, $Y=\hat{Y}/F$ is a $\boldsymbol{%
\Sigma }_{1}^{1}$-definable set such that $F$ is Borel, and $f:X\rightarrow Y
$ is a Borel-definable bijection, then $f$ is an isomorphism in the category
of $\boldsymbol{\Sigma }_{1}^{1}$-definable sets; see also \cite[Proposition
3.7]{bergfalk_definable_2024-1}. Furthermore, by \cite[Proposition 3.7]%
{bergfalk_definable_2024-1}, \emph{Borel-definable sets }form a \emph{%
replete }full subcategory of the category of $\boldsymbol{\Sigma }_{1}^{1}$%
-definable sets.

\subsection{Abelian groups with a Polish cover}

A \emph{Polish group }is a topological group whose topology is Polish.
Abelian Polish groups form a countably complete quasi-abelian category $%
\mathbf{PAb}$, whose morphisms are the continuous group homomorphisms, and
whose strict morphisms are the continuous group homomorphisms with closed
image. Its left heart $\mathrm{LH}\left( \mathbf{PAb}\right) $ is
equivalently described in \cite{lupini_looking_2024} as the category of 
\emph{abelian groups with a Polish cover }considered in \cite%
{bergfalk_definable_2024,bergfalk_definable_2024-1} and also by Moore in 
\cite{moore_group_1976} (where they are called \emph{pseudo-polonais abelian
groups}).

An abelian group with a Polish cover is an abelian group $G$ explicitly
presented as a quotient $\hat{G}/N$ where $\hat{G}$ is an abelian Polish
group (the \textquotedblleft Polish cover\textquotedblright ) and $N$ is a 
\emph{Polish subsgroup} of $\hat{G}$. This means that $N$ is itself a Polish
group with respect to a topology that makes the inclusion $N\rightarrow \hat{%
G}$ continuous or, equivalently, $N$ is the image of a continuous group
homomorphism $H\rightarrow \hat{G}$ from some Polish group $H$. In
particular, this implies that $N$ is a Borel subset of $\hat{G}$. Abelian
groups with a Polish cover are the objects of a countably complete category,
whose morphisms $G\rightarrow H$ are the group homomorphism that are \emph{%
Borel-definable}, i.e., induced by a Borel function (which is not
necessarily a group homomorphism) between the Polish covers. It is proved in 
\cite{lupini_looking_2024} that this is an abelian category, which is in
fact equivalent to the left heart $\mathrm{LH}\left( \mathbf{PAb}\right) $
of $\mathbf{PAb}$. The canonical inclusion $\mathbf{PAb}\rightarrow \mathrm{%
LH}\left( \mathbf{PAb}\right) $ is obtained by identifying an abelian Polish
group $G$ with the group with a Polish cover $G/N$ where $N$ is the trivial
subgroup of $G$.

\subsection{Homogeneous spaces with a Polish cover}

Abelian groups with a Polish cover are a particular instance of \emph{%
homogeneous spaces with a Polish cover}. By definition, this is a
homogeneous space $X=G/N$ where $G$ is a Polish group and $N\subseteq G$ is
a Polish subgroup. We regard this as a pointed space, with the distinguished
point $\ast \in G/N$ corresponding to $N$. When $N$ is normal in $G$, $X$ is
called a group with a Polish cover. A subspace with a Polish cover of $X$ is
a subspace $Y=H/N$ for some Polish subgroup of $G$ containing $N$. Defining
the Borel class and Borel rank of $Y$ in $X$ to be the Borel class and Borel
rank of $H$ in $G$ yields a canonical notion of complexity.

It was proved by Solecki in \cite{solecki_polish_1999} that a homogeneous
space with a Polish cover $X$ has a canonical chain of Polish subgroups
indexed by countable ordinals, which we denote by $\mathrm{Ph}^{\alpha
}\left( X\right) $ for $\alpha <\omega _{1}$ and call \emph{phantom
subspaces }(or Solecki subspaces). These were characterized in terms of
complexity in \cite{farah_borel_2006}, where it is shown that $\mathrm{Ph}%
^{\alpha }\left( X\right) $ is the smallest $\boldsymbol{\Pi }_{1+\alpha
+1}^{0}$ pointed subspace of $X$; see also \cite[Section 4.3]%
{lupini_looking_2024}. In particular, $\mathrm{Ph}^{0}\left( X\right) $ is
the smallest closed pointed subspace of $X$ (with respect to the quotient
topology). By convention, we set $\mathrm{Ph}^{-1}\left( X\right) =X$. When $%
X$ is a group with a Polish cover, $\mathrm{Ph}^{\alpha }\left( X\right) $
is a subgroup with a Polish cover of $X$ for every $\alpha <\omega _{1}$.

\begin{definition}
\label{Definition:phantom-length-bound}Let $X$ be a homogeneous space with a
Polish cover, and $\alpha <\omega _{1}$ and $\lambda <\omega _{1}$ zero or
limit. Then:

\begin{itemize}
\item $X$ has phantom length at most $\alpha $ if $\mathrm{Ph}^{\alpha
}X=\left\{ \ast \right\} $

\item $X$ has phantom length $\alpha +1+1/2$ if $\mathrm{Ph}^{\alpha
+2}X=\left\{ \ast \right\} $ and $\left\{ \ast \right\} \in D(\boldsymbol{%
\Pi }_{2}^{0})\left( \mathrm{Ph}^{\alpha }X\right) $;

\item $X$ has phantom length at most $\lambda +1/2$ if $\mathrm{Ph}^{\lambda
+1}X=\left\{ \ast \right\} $ and $\left\{ \ast \right\} \in \boldsymbol{%
\Sigma }_{2}^{0}(\mathrm{Ph}^{\lambda }X)$;

\item $X$ has phantom length at most $\lambda +1/2+\varepsilon $ if $\mathrm{%
Ph}^{\lambda +1}X=\left\{ \ast \right\} $ and $\left\{ \ast \right\} \in D(%
\boldsymbol{\Pi }_{2}^{0})(\mathrm{Ph}^{\lambda }X)$.
\end{itemize}

Then \emph{phantom length }of $X$ is the least $\alpha \in \omega _{1}^{%
\mathbf{Pol}}$ such that $X$ has phantom length at most $\alpha $.
\end{definition}

By \cite[Theorem 6.1]{lupini_complexity_2025}, the phantom length of $X$
completely determines the complexity class of $\left\{ \ast \right\} $ in $X$%
, and vice versa. Furthermore, basepoint-preserving Borel-definable
bijections between homogeneous spaces with a Polish cover preserve the
phantom length, with the exception of the values $1/2$ and $1/2+\varepsilon $%
, which can be interchanged \cite[Theorem 3.3]{lupini_complexity_2025}. If $%
X=G/N$ where $N$ is \emph{non-Archimedean}, then the phantom length of $X$
is in $\omega \lbrack 1/2]$ by \cite[Theorem 6.1]{lupini_complexity_2025};
see Section \ref{Subsection:subcategories}.

\subsection{Subcategories\label{Subsection:subcategories}}

Recall that a Polish group $H$ is \emph{non-Archimedean} when it has a basis
of identity neighborhoods consisting of subgroups. When $H$ is abelian, this
is equivalent to being\emph{\ pro-countable}.

\begin{definition}
A Polish group is \emph{pro-lc} if it is the limit of an \emph{epimorphic}
tower of \emph{locally compact} Polish groups.
\end{definition}

Thus, a Polish group $G$ is pro-lc if for every identity neighborhood $U$ of 
$G$ there exists a closed normal subgroup $H$ of $G$ contained in $U$ such
that $G/H$ is locally compact. An \emph{abelian }Polish group is pro-lc if
and only if it is pro-Lie \cite{hofmann_lie_2007}, i.e., the limit of a
tower of (abelian) Lie groups \cite{casarosa_homological_2026}. An abelian
Lie group is necessarily of the form $V\oplus T\oplus C$ where $V$ is a
finite-dimensional vector group, $T$ is a finite-dimensional torus group,
and $C$ is a countable abelian group.

A strictly full subcategory $\mathcal{B}$ of a quasi-abelian category $%
\mathcal{A}$ is a quasi-abelian subcategory if it is quasi-abelian and the
inclusion $\mathcal{B}\rightarrow \mathcal{A}$ preserves finite limits and
finite colimits \cite[Definition 2.7]{casarosa_homological_2026}. When
furthermore $\mathcal{B}$ is closed under extensions, then $\mathcal{B}$ is
a thick subcategory of $\mathcal{A}$ \cite[Definition 2.7]%
{casarosa_homological_2026}; see also \cite[Definition 6.15]%
{lupini_looking_2024}. It is shown in \cite[Theorem 6.17]%
{lupini_looking_2024} and \cite[Theorem 3.17]{casarosa_homological_2026}
that both pro-countable abelian Polish groups and pro-Lie abelian Polish
groups form thick subcategories of $\mathbf{PAb}$.

If $\mathcal{B}$ is a thick subcategory of $\mathbf{PAb}$, then one defines
a group with a cover in $\mathcal{B}$ to be a group with a Polish cover $G=%
\hat{G}/N$ where \emph{both} $N$ and $\hat{G}$ are in $\mathcal{B}$ \cite[%
Section 6.3]{lupini_looking_2024}. In this case, one has that $\mathrm{LH}%
\left( \mathcal{B}\right) $ is equivalent to the full\emph{\ }subcategory of 
$\mathrm{LH}\left( \mathbf{PAb}\right) $ spanned by the groups with a cover
in $\mathcal{B}$ \cite[Proposition 6.16]{lupini_looking_2024}. As $\mathcal{B%
}$ is thick in $\mathbf{PAb}$, $\mathrm{LH}\left( \mathcal{B}\right) $ is an
abelian subcategory of $\mathrm{LH}\left( \mathbf{PAb}\right) $.

\subsection{Complexity}

We recall the following result from \cite[Theorem 3.5.3]{gao_invariant_2009}.

\begin{lemma}
\label{Lemma:Mackey}Let $G$ be a closed subgroup of a Polish group $H$.
Suppose that $X$ is a standard Borel $G$-space. Then there exists a Borel $H$%
-space $Y$ containing $X$ as a Borel subset such that:

\begin{enumerate}
\item the inclusion $X\rightarrow Y$ is $G$-equivariant;

\item every $H$-orbit of $Y$ contains exactly one $G$-orbit of $X$.
\end{enumerate}
\end{lemma}

Let $S_{\infty }$ be the Polish group of permutations of $\mathbb{N}$
endowed with the topology of pointwise convergence. Recall \cite%
{jackson_countable_2002,hjorth_classification_2000} that a Borel equivalence
relation $E$ on a standard Borel space $X$ is:

\begin{itemize}
\item \emph{countable }if it is the orbit equivalence relation of a standard
Borel $G$-space for some \emph{countable }group $G$;

\item \emph{hyperfinite }if it is the orbit equivalence relation of a
standard Borel $G$-space where $G$ is the additive group of $\mathbb{Z}$;

\item \emph{essentially hyperfinite }if it is Borel-reducible to a Borel
equivalence relation that is hyperfinite;

\item \emph{essentially countable }if it is Borel-reducible to a Borel
equivalence relation that is hyperfinite;

\item \emph{classifiable by countable structures }if it is Borel-reducible
to the orbit equivalence relation associated with a standard Borel $%
S_{\infty }$-space.
\end{itemize}

One denotes by $E_{0}$ the relation on $2^{\mathbb{N}}$ of tail equivalence
of binary sequences. Its countable product $E_{0}^{\mathbb{N}}$ is the
equivalence relation on $\left( 2^{\mathbb{N}}\right) ^{\mathbb{N}}$ defined
by%
\begin{equation*}
\left( \boldsymbol{x}_{i}\right) E_{0}^{\mathbb{N}}\left( \boldsymbol{y}%
_{i}\right) \Leftrightarrow \forall i\in \mathbb{N}\text{, }\boldsymbol{x}%
_{i}E_{0}\boldsymbol{y}_{i}\text{.}
\end{equation*}%
If $G$ is a non-Archimedean Polish space, then $G$ is isomorphic to a closed
subgroup of $S_{\infty }$. Thus, if $X$ is a standard Borel $G$-space, then $%
E_{G}^{X}$ is classifiable by countable structures by Lemma \ref%
{Lemma:Mackey}. Furthermore, if $E_{G}^{X}$ is potentially $\boldsymbol{%
\Sigma }_{3}^{0}$, then it is potentially $\boldsymbol{\Sigma }_{2}^{0}$ 
\cite[Theorem 4.1]{hjorth_borel_1998}.\ If $E_{G}^{X}$ is potentially $%
\boldsymbol{\Pi }_{2}^{0}$, then it is potentially $\boldsymbol{\Pi }_{1}^{0}
$ \cite[Theorem 3.1]{hjorth_borel_1998}.

By Lemma \ref{Lemma:Mackey} together with \cite[Corollary 4.5 and Theorem 4.6%
]{allison_countable_2025} one obtains the following:

\begin{lemma}[Allison]
\label{Lemma:allison}Let $G$ be a pro-lc Polish group, and $X$ is a Polish $%
G $-space. Let $E_{G}^{X}$ be the corresponding orbit equivalence relation.
Then:

\begin{enumerate}
\item if $E_{G}^{X}$ is $\boldsymbol{\Pi }_{3}^{0}$, then $E$ is Borel
reducible to $E_{0}^{\mathbb{N}}$ and in particular classifiable by
countable structures;

\item if $E$ is a countable Borel equivalence relation such that $E$ is
Borel reducible $E_{G}^{X}$, then $E$ is Borel reducible to $E_{0}$.
\end{enumerate}
\end{lemma}

The following result is the particular case of \cite[Corollary 2.8]%
{allison_non-archimedean_2020} when $\Gamma $ is $\boldsymbol{\Sigma }%
_{2}^{0}$.

\begin{lemma}[Allison]
\label{Lemma:Allison2}Let $G$, $H$ be Polish groups, $X$ be a Polish $G$%
-space, and $Y$ be a Polish $H$-space. Let also $E$ be an equivalence
relation on a standard Borel space. Suppose that:

\begin{enumerate}
\item $E$ is Borel-reducible to $E_{G}^{X}$;

\item $E$ is Borel-reducible to $E_{H}^{Y}$;

\item $E_{H}^{Y}$ is potentially $\boldsymbol{\Sigma }_{2}^{0}$.
\end{enumerate}

Then there exists a Polish $G$-space $Z$ such that:

\begin{itemize}
\item $E_{G}^{Z}$ is $\boldsymbol{\Sigma }_{2}^{0}$, and

\item $E$ is Borel-reducible to $E_{G}^{Z}$.
\end{itemize}
\end{lemma}

\begin{corollary}
\label{Corollary:Allison}Let $E$ be an equivalence relation on a standard
Borel space. Let also $H$ be a Polish group, and $Y$ be a Polish $H$-space.
Suppose that:

\begin{enumerate}
\item $E$ is classifiable by countable structure;

\item $E$ is the orbit equivalence relation associated with an action of a
pro-lc Polish group on a Polish space;

\item $E$ is Borel-reducible to $E_{H}^{Y}$;

\item $E_{H}^{Y}$ is potentially $\boldsymbol{\Sigma }_{3}^{0}$.
\end{enumerate}

Then $E$ is essentially hyperfinite.
\end{corollary}

\begin{proof}
Let $G=S_{\infty }$ be the Polish groups of permutations of $\mathbb{N}$.
Since $E$ is classifiable by countable structures, it is Borel-reducible to $%
E_{G}^{X}$ for some Polish $G$-space $X$. By Lemma \ref{Lemma:Allison2}
there exists a Polish $G$-space $Z$ such that $E_{G}^{Z}$ is $\boldsymbol{%
\Sigma }_{3}^{0}$ and $E$ is Borel-reducible to $E_{G}^{Z}$. By \cite[%
Theorem 12.5.7]{gao_invariant_2009}, $E_{G}^{Z}$ is essentially countable.
By \cite[Lemma 3.9]{kechris_borel_2016}, there exist a countable Borel
equivalence relation $F$ such that $E$ is Borel-reducible to $F$ and $F$ is
Borel-reducible to $E$. Thus, by Lemma \ref{Lemma:allison}(2), $F$ is
essentially hyperfinite. Thus, the same holds for $E$.
\end{proof}

For a homogeneous space with a Polish cover $X=G/N$, we let $=_{G/N}$ be the
coset relation of $G$ within $N$. Then $=_{G/N}$ is smooth if and only if $N$
is closed in $G$; see \cite[Proposition 4.12]{lupini_looking_2024}.

\begin{proposition}
\label{Proposition:Pi3}Suppose that $G$ is a Polish group, and $N\subseteq G$
is a pro-lc Polishable subgroup. Then the following assertions are
equivalent:

\begin{enumerate}
\item $=_{G/N}$ is potentially $\boldsymbol{\Pi }_{3}^{0}$;

\item $=_{G/N}$ is Borel reducible to $E_{0}^{\mathbb{N}}$;

\item $N$ is $\boldsymbol{\Pi }_{3}^{0}$ in $G$.
\end{enumerate}
\end{proposition}

\begin{proof}
(1)$\Rightarrow $(2) follows from Lemma \ref{Lemma:allison}; (2)$\Rightarrow 
$(3) follows from the fact that $E_{0}^{\omega }$ is $\boldsymbol{\Pi }%
_{3}^{0}$; (1)$\Rightarrow $(3) follows from \cite[Theorem 3.3(3)]%
{lupini_complexity_2025}.
\end{proof}

Let $G$ be a Polish group, and $N\subseteq G$ be a Polishable subgroup.
Recall that if $N$ is $\boldsymbol{\Pi }_{3}^{0}$ in $G$, then $N$ has a
basis of neighborhoods of the identity $V$ satisfying $\overline{V}^{G}\cap
N=V$; see \cite[Section 4]{lupini_complexity_2025}. It follows from this and
the Baire Category Theorem for $N$ that if $N$ is $\boldsymbol{\Sigma }%
_{2}^{0}$ in $G$, then $N$ has a basis of neighborhoods of the identity that
are closed in $G$. For $A,B\subseteq G$ set%
\begin{equation*}
AB:=\left\{ xy:x\in A\text{ and }y\in B\right\} \subseteq G\text{.}
\end{equation*}%
The proof of the implication (4)$\Rightarrow $(2) in the following lemma was
suggested by Allison.

\begin{proposition}
\label{Proposition:Sigma2}Suppose that $G$ is a Polish group, and $N$ is a
pro-lc Polishable subgroup of $G$. Then the following assertions are
equivalent:

\begin{enumerate}
\item $=_{G/N}$ is potentially $\boldsymbol{\Sigma }_{2}^{0}$;

\item $=_{G/N}$ is Borel reducible to $E_{0}$;

\item $N$ is $\boldsymbol{\Sigma }_{2}^{0}$ in $G$;

\item $N$ is $\boldsymbol{\Sigma }_{3}^{0}$ in $G$.
\end{enumerate}
\end{proposition}

\begin{proof}
The equivalence of (1) and (4) is the content of \cite[Theorem 3.3(2)]%
{lupini_complexity_2025}.

(4)$\Rightarrow $(2) By \cite[Theorem 6.1]{lupini_complexity_2025}, $N$ is $%
D(\boldsymbol{\Pi }_{2}^{0})$ in $G$. In particular $=_{G}$ is $\boldsymbol{%
\Pi }_{3}^{0}$. Thus, by Lemma \ref{Lemma:allison}, $=_{G}$ is Borel
reducible to $E_{0}^{\mathbb{N}}$ and, in particular, classifiable by
countable structures. The conclusion thus follows from Corollary \ref%
{Corollary:Allison}.

(4)$\Rightarrow $(3) Assume initially that $N$ is the product of locally
compact Polish abelian groups. As in the proof of \cite[Proposition 6.9]%
{lupini_complexity_2025}, the proof of \cite[Lemma 6.4]%
{lupini_complexity_2025} shows that there exists a neighborhood $U$ of the
identity in $N$ satisfying $U=U^{-1}$ that is $\boldsymbol{\Pi }_{2}^{0}$ in 
$G$. We can write $N$ as the (internal) direct product $N=H\times N_{0}$
where $H\subseteq N$ is locally compact and $N_{0}\subseteq U\subseteq N$.
Since $H$ is locally compact, it is $\boldsymbol{\Sigma }_{2}^{0}$ in $G$.
Let us consider the Vaught transform \cite{gao_invariant_2009} with respect
to the action $N\curvearrowright G$ by translation. Then $M:=U^{\ast N}$ is $%
\boldsymbol{\Pi }_{2}^{0}$ in $G$ \cite[Theorem 3.2.7]{gao_invariant_2009}
and $N$-invariant \cite[Proposition 3.2.3]{gao_invariant_2009}. In
particular, we have that $U^{\ast N}$ is a subsemigroup of $N$. Since $%
U=U^{-1}$, it follows that $M=M^{-1}$ and $M$ is a subgroup of $G$. This
implies that $U^{\ast N}$ is a \emph{closed }subgroup of $G$, which contains 
$N_{0}$. If $K$ is a compact neighborhood of the identity of $H$, then $KM$
is closed in $G$. Since $H$ is a countable union of translates of $K$, it
follows that $N$ is a countable union of translates of $KM$, and hence $%
\boldsymbol{\Sigma }_{2}^{0}$ in $G$.

Suppose now that $N$ is an arbitrary pro-lc locally compact group. Then
there exists a product of locally compact Polish abelian groups $N^{\prime }$
that contains $N$ as a closed subgroup. Consider now the quotient $G^{\prime
}$ of $N^{\prime }\times G$ by the closed subgroup 
\begin{equation*}
\left\{ \left( x,y\right) \in N^{\prime }\times G:x=y\in N\right\} \text{.}
\end{equation*}%
Then we can canonically identify $N^{\prime }$ as a $\boldsymbol{\Sigma }%
_{3}^{0}$ Polishable subgroup of $G^{\prime }$ and $G$ as a closed subgroup
of $G^{\prime }$. Since $N^{\prime }$ is a product of locally compact Polish
abelian groups, we conclude by the above applied to $N^{\prime }\subseteq
G^{\prime }$ that $N^{\prime }$ is in fact $\boldsymbol{\Sigma }_{2}^{0}$ in 
$G^{\prime }$. Since $G$ is closed in $G^{\prime }$, its topology is the
subspace topology inherited from $G^{\prime }$. Therefore, $N=N^{\prime
}\cap G$ is $\boldsymbol{\Sigma }_{2}^{0}$ in $G$.
\end{proof}

\begin{corollary}
The set $\omega \lbrack 1/2]$ provides a complete list of the possible
phantom lengths of homogeneous spaces with a Polish cover $X=G/N$ where $N$
is pro-lc.
\end{corollary}

\begin{proof}
This follows from Proposition \ref{Proposition:Sigma2} and \cite[Theorem 6.1]%
{lupini_complexity_2025}.
\end{proof}

\begin{corollary}
The phantom length of a homogeneous space with a Polish cover $X=G/N$ where $%
N$ is pro-lc is invariant under Borel-definable basepoint-preserving
bijections.
\end{corollary}

\begin{proof}
This follows from Proposition \ref{Proposition:Sigma2} and \cite[Theorem 3.3]%
{lupini_complexity_2025}.
\end{proof}

\subsection{Hereditarily countable sets}

For every $\alpha <\omega _{1}$, a canonical Polish $S_{\infty }$-space $%
X_{\alpha }$ is constructed in \cite{hjorth_borel_1998} such that the
potential complexity class of the corresponding orbit equivalence relation $%
E_{S_{\infty }}^{X_{\alpha }}$ is $\Gamma _{\alpha }$. Furthermore, for $%
\alpha \geq 2$ there exists a Borel $S_{\infty }$-invariant subspace $%
X_{\alpha -1/2}$ such that the potential complexity class of $E_{S_{\infty
}}^{X_{\alpha }}$ is $\Gamma _{\alpha -1/2}$.

In fact, $E_{S_{\infty }}^{X_{\alpha }}$ can be seen as the relation of 
\emph{isomorphism }for a class of countable structures in a first-order
language providing a canonical parametrization of \emph{hereditarily
countable sets of rank} at most $1+\alpha $. We let%
\begin{equation*}
\wp ^{\alpha }\left( \mathbb{N}\right) :=X_{\alpha }/S_{\infty }
\end{equation*}%
be the corresponding $\boldsymbol{\Sigma }_{1}^{1}$-definable set. Set also%
\begin{equation*}
\wp ^{\alpha -1/2}\left( \mathbb{N}\right) :=X_{\alpha -1/2}/S_{\infty }
\end{equation*}%
for $\alpha \geq 2$, which we see as a parametrization of hereditarily
countable sets of rank at most $\alpha -1/2$. For $\alpha =1/2$ define%
\begin{equation*}
\wp ^{1/2}\left( \mathbb{N}\right) :=X/E
\end{equation*}%
where $E$ is the countable Borel equivalence relation on a Polish space $X$
of maximum Borel complexity \cite{jackson_countable_2002}. Let also $\wp
^{0}\left( \mathbb{N}\right) $ be the Polish space of subsets of $\mathbb{N}$%
.

\begin{definition}
\label{Definition:hereditarily-countable}Let $X=\hat{X}/E$ be a $\boldsymbol{%
\Sigma }_{1}^{1}$-definable set, and $\alpha \in \omega _{1}[1/2]$. The
elements of $\hat{X}$ are \emph{classifiable up to} $E$ \emph{by
hereditarily countable sets of rank} $1+\alpha $, and $X$ can be \emph{%
parametrized by hereditarily countable sets of rank} $1+\alpha $, if there
exists a Borel-definable injection $X\rightarrow \wp ^{\alpha }\left( 
\mathbb{N}\right) $.
\end{definition}

The following result is obtained in \cite[Theorem 2 and Theorem 3]%
{hjorth_borel_1998}; see also \cite[Lemma 13.19 and Theorem 13.21]%
{casarosa_phantom_2025}.

\begin{theorem}
\label{Theorem:hereditarily-countable}Let $E$ be a $\boldsymbol{\Sigma }%
_{1}^{1}$ equivalence relation on a Polish space $\hat{X}$ that is
classifiable by countable structures. Then for every $\alpha \in \omega
_{1}[1/2]$, $E$ has potential complexity $\Gamma _{\alpha }$ if and only if
the points of $\hat{X}$ are classifiable up to $E$ by hereditarily countable
sets of rank $1+\alpha $.
\end{theorem}

\section{Phantom maps\label{Section:phantom-maps}}

In this section we recall the notion of \emph{phantom map }(of the second
kind) \cite%
{mcgibbon_phantom_1995,mcgibbon_phantom_1994,mcgibbon_phantom_1997,gray_universal_1993,mcgibbon_numerical_2001}
and its \emph{higher order} versions \cite{ha_higher_2003}; see also \cite%
{lee_phantom_2004,ha_gray_2005,tsakanikas_rationalization_2010,iriye_gray_2010}%
. We then obtain a \emph{bound} on the order of nontrivial phantom maps.

\subsection{Homotopy}

Let $X,Y$ be locally compact Polish spaces. Suppose that $f,g:X\rightarrow Y$
are continuous maps. Then $f,g$ are \emph{homotopic} if there exists a
continuous function $h:X\times \left[ 0,1\right] \rightarrow Y$ such that $%
f=h\left( -,0\right) $ and $g=h\left( -,1\right) $. This defines an
equivalence relation $\sim $ on the space $\mathbf{LC}\left( X,Y\right) $ of
continuous maps $X\rightarrow Y$, which is a Polish space when endowed with
the compact-open topology. The corresponding quotient space is denoted by $%
[X,Y]$.

It is shown in \cite[Theorem 4.15]{bergfalk_definable_2024-1} that when $X\ $%
is a locally compact Polish space and $P$ is a homotopy polyhedron, the
relation of homotopy for maps $X\rightarrow P$ is an \emph{analytic} and 
\emph{idealistic }equivalence relation. When furthermore $P$ is a \emph{%
polyhedral }$H$-\emph{group }\cite[Section 5.1]{bergfalk_definable_2024-1},
the relation of homotopy for maps $X\rightarrow P$ is Borel \cite[Theorem
7.13]{bergfalk_definable_2024-1}, and $[X,P]$ is a Borel-definable set. The $%
H$-group operations on $P$ endows $[X,P]$ with a group structure, turning it
into a group object in the category of Borel-definable sets (\emph{%
Borel-definable group}).

\subsection{Phantom maps}

Let $X$ be a locally compact Polish space, and $\left( P,\ast \right) $ a
pointed homotopy polyhedron. A \emph{phantom map }(of the second kind) from $%
X$ to $P$ is a continuous function $f:X\rightarrow P$ such that, for every
compact subset $C$ of $X$, $f|_{C}$ is \emph{nullhomotopic}, i.e., homotopic
to the function $\ast :C\rightarrow P$ constantly equal to $\ast $ \cite%
{mcgibbon_phantom_1995}.

\emph{Higher order }phantom maps have been introduced in \cite%
{ha_higher_2003}. These are defined by recursion on countable ordinals as
follows. A map $f:X\rightarrow P$ is phantom of order at least $0$ if and
only if it is a phantom map in the sense defined above, i.e., for every
compact $C\subseteq X$, the restriction $f|_{C}$ is homotopic to the
function $\ast :C\rightarrow P$ constantly equal to $\ast $.\ This is
equivalent to the assertion that there exists a map $g:X/C\rightarrow P$
such that the diagram%
\begin{equation*}
\begin{array}{ccc}
X & \rightarrow & P \\ 
\downarrow & \nearrow &  \\ 
X/C &  & 
\end{array}%
\end{equation*}%
commutes \emph{up to homotopy}. Here, $X/C$ is the quotient of $X$ obtained
by identifying $C$ to a point, and $X\rightarrow X/C$ is the canonical
quotient map.

Recursively, for every countable ordinal $\alpha \geq 1$, one lets $f$ be
phantom of order at least $\alpha $ if for every compact $C\subseteq X$ and
for every $\beta <\alpha $ there exists a \emph{phantom map }$%
g:X/C\rightarrow P$ \emph{of order at least} $\beta $ such that the diagram%
\begin{equation*}
\begin{array}{ccc}
X & \overset{f}{\rightarrow } & P \\ 
\downarrow  & \nearrow  &  \\ 
X/C &  & 
\end{array}%
\end{equation*}%
commutes \emph{up to homotopy}. We denote by $\mathrm{Ph}^{\alpha }[X,P]$
the space of maps $X\rightarrow P$ that are phantom of order at least $%
\alpha \in \omega _{1}$. We also set%
\begin{equation*}
\mathrm{Ph}^{-1}[X,P]:=[X,P]\text{.}
\end{equation*}%
Notice that, by definition, for a limit ordinal $\lambda $, 
\begin{equation*}
\mathrm{Ph}^{\lambda }[X,P]=\bigcap_{n\in \omega }\mathrm{Ph}^{\lambda
_{n}}[X,P]\text{.}
\end{equation*}%
We can write for every countable ordinal $\alpha \geq 1$: 
\begin{equation}
\mathrm{Ph}^{\alpha }[X,P]=\bigcap_{\beta <\alpha }\bigcap_{C\subseteq X}%
\mathrm{Ran}\left( \mathrm{Ph}^{\beta }[X/C,P]\rightarrow \lbrack
X,P]\right) \text{.\label{Equation:phantom}}
\end{equation}%
where $C$ ranges among the compact subspaces of $X$. Here,%
\begin{equation*}
\lbrack X/C,P]\rightarrow \lbrack X,P]
\end{equation*}%
is the function induced by the quotient map $X\rightarrow X/C$. Adopting the
notation from Section \ref{Subsection:notation}, we can write the following
expression, subsuming Equation \ref{Equation:phantom} for $\alpha \geq 1$.
For every $\alpha <\omega _{1}$:%
\begin{equation*}
\mathrm{Ph}^{\alpha }[X,P]=\bigcap_{n\in \omega }\bigcap_{C\subseteq X}%
\mathrm{Ran}\left( \mathrm{Ph}^{\left( \alpha -1\right)
_{n}}[X/C,P]\rightarrow \lbrack X,P]\right) 
\end{equation*}

\subsection{Derived towers\label{Section:derived}}

Let $\boldsymbol{G}:=\left( G^{\left( n\right) }\right) _{n\in \omega }$ be
a tower of (not necessarily abelian) countable groups. Consider the Polish
group $G:=\prod_{n\in \omega }G^{\left( n\right) }$. Then $G$ admits a
continuous action on itself by setting 
\begin{equation*}
\left( g_{n}\right) \cdot \left( x_{n}\right) =\left( g_{n}\cdot x_{n}\cdot
p^{\left( n,n+1\right) }\left( g_{n+1}\right) \right) _{n\in \omega }\text{.}
\end{equation*}%
This defines a pointed $\mathbf{\Sigma }_{1}^{1}$ definable set (the
distinguished element corresponding to the identity element of $G$), which
is denoted by $\mathrm{\mathrm{\mathrm{li}}m}_{n}^{1}G^{\left( n\right) }$
or $\mathrm{lim}^{1}\boldsymbol{G}$. The assignment $\left( G_{n}\right)
\mapsto \mathrm{\mathrm{lim}}_{n}^{1}G_{n}$ defines a functor from the
category of towers of countable groups to the category of pointed $\mathbf{%
\Sigma }_{1}^{1}$-definable sets. When $\boldsymbol{G}$ is a tower of
(additively denoted) countable \emph{abelian }groups, the orbit equivalence
relation of the action $G\curvearrowright G$ defining $\mathrm{\mathrm{lim}}%
^{1}\boldsymbol{G}$ is the coset relation with respect to the Polishable
subgroup $\left\{ \left( g_{n}-p^{\left( n,n+1\right) }\left( g_{n+1}\right)
\right) :\left( g_{n}\right) \in G\right\} $. Thus, we have that $\mathrm{%
\mathrm{lim}}^{1}\boldsymbol{G}$ is an abelian group with a Polish cover. In
the abelian case, $\mathrm{lim}^{1}$ describes the derived functor of the
limit functor $\mathrm{lim}$ for towers \cite[Chapter 11]%
{mardesic_strong_2000}.

The \emph{derived tower }$\boldsymbol{G}^{\prime }$ of $\boldsymbol{G}$ is
the tower $\left( G^{\prime \left( n\right) }\right) $ defined by setting%
\begin{equation*}
G^{\prime \left( n\right) }:=\bigcap_{k>n}p^{\left( n,k\right) }(G^{\left(
k\right) })
\end{equation*}%
for $n\in \omega $. One then defines by recursion on $\alpha <\omega _{1}$
the $\alpha $-th derived tower $\boldsymbol{G}_{\alpha }$ of $\boldsymbol{G}$
by setting $\boldsymbol{G}_{0}:=\boldsymbol{G}$, $\boldsymbol{G}_{\alpha
+1}:=\boldsymbol{G}_{\alpha }^{\prime }$ for $\alpha <\omega _{1}$, and for
a limit ordinal $\lambda $, $\boldsymbol{G}_{\lambda }=(G_{\lambda
_{n}}^{\left( n\right) })$. Notice that%
\begin{equation*}
G_{\lambda +1}^{\left( n\right) }=\bigcap_{\beta <\lambda }G_{\beta
}^{\left( n\right) }\text{.}
\end{equation*}%
The towers $\boldsymbol{G}_{\alpha }$ are also called \emph{image towers} in 
\cite{boardman_conditionally_1999,ha_higher_2003}, where a slightly
different indexing is used. It is easy to verify by induction on $\alpha
<\omega _{1}$ that $\boldsymbol{G}\mapsto \boldsymbol{G}_{\alpha }$ is a
subfunctor of the identity on the category of towers of countable groups. In
particular, if $\boldsymbol{G}$ and $\boldsymbol{H}$ are isomorphic towers,
then $\boldsymbol{G}_{\alpha }$ and $\boldsymbol{H}_{\alpha }$ are
isomorphic towers for every $\alpha <\omega _{1}$.

Notice that the inclusion $\boldsymbol{G}^{\prime }\rightarrow \boldsymbol{G}
$ induces an injective basepoint-preserving definable function $\mathrm{%
\mathrm{\mathrm{li}}m}^{1}\boldsymbol{G}^{\prime }\rightarrow \mathrm{%
\mathrm{lim}}^{1}\boldsymbol{G}$. It follows by induction on $\alpha <\omega
_{1}$ that the inclusion $\boldsymbol{G}_{\alpha }\rightarrow \boldsymbol{G}$
induces an injective basepoint-preserving definable function $\mathrm{%
\mathrm{lim}}^{1}\boldsymbol{G}_{\alpha }\rightarrow \mathrm{\mathrm{lim}}%
^{1}\boldsymbol{G}$, which allows one to regard $\mathrm{\mathrm{lim}}^{1}%
\boldsymbol{G}_{\alpha }$ as a pointed subset of\textrm{\ lim}$^{1}%
\boldsymbol{G}$.

A tower of countable groups $\boldsymbol{G}$:

\begin{itemize}
\item is epimorphic if $p^{\left( n,n+1\right) }:G^{\left( n+1\right)
}\rightarrow G^{\left( n\right) }$ is surjective for every $n\in \mathbb{N}$;

\item is monomorphic if $p^{\left( n,n+1\right) }:G^{\left( n+1\right)
}\rightarrow G^{\left( n\right) }$ is injective for every $n\in \mathbb{N}$;

\item is essentially monomorphic if it is isomorphic to an injective tower;

\item is essentially epimorphic, or satisfies the \emph{Mittag-Leffler
property }if for every $n\in \omega $ the sequence $\left( p^{\left(
n,k\right) }\left( G_{k}\right) \right) _{k>n}$ of subgroups of $G_{n}$ is
eventually constant.
\end{itemize}

A tower of Polish groups $\boldsymbol{G}$ satisfies the Mittag-Leffler
property if and only if $\boldsymbol{G}$ is isomorphic to a tower $%
\boldsymbol{H}$ with $\boldsymbol{H}^{\prime }=\boldsymbol{H}$ \cite[Chapter
II, Section 6.2, Theorem 7]{mardesic_shape_1982}, which implies $\mathrm{%
\mathrm{lim}}^{1}\boldsymbol{G}=0$. Conversely, if $\boldsymbol{G}$ is a
tower of countable groups such that $\mathrm{lim}^{1}\boldsymbol{G}=\left\{
\ast \right\} $, then it satisfiest he Mittag-Leffler property \cite[Chapter
II, Section 6.2, Theorems 10 and 11]{mardesic_shape_1982}. 

When $\boldsymbol{G}$ is a tower of countable groups such that $\boldsymbol{G%
}^{\prime }$ satisfies the Mittag--Leffler property, the canonical
surjective Borel-definable pointed map 
\begin{equation*}
\mathrm{lim}^{1}\boldsymbol{G}\rightarrow \mathrm{lim}_{n}\left( \mathrm{lim}%
_{k}G_{n}/p^{\left( n,k\right) }G_{k}\right) 
\end{equation*}%
is injective; see \cite[Section 7.2]{casarosa_phantom_2025}. Thus, in this
case $\mathrm{lim}^{1}\boldsymbol{G}$ is $\boldsymbol{\Pi }_{3}^{0}$%
-definable. Furthermore, if $\boldsymbol{G}$ is a tower of countable groups,
then $\mathrm{lim}^{1}\boldsymbol{G}$ is $\boldsymbol{\Sigma }_{2}^{0}$%
-definable if and only if $\boldsymbol{G}$ is essentially monomorphic \cite[%
Lemma 7.5]{casarosa_phantom_2025}.

\begin{lemma}
\label{Lemma:countable-rank}Let $\boldsymbol{A}=\left( A^{\left( n\right)
}\right) $ be a tower of Polish groups. Then there exists $\alpha <\omega
_{1}$ such that $\boldsymbol{A}_{\alpha }=\boldsymbol{A}_{\alpha +1}$, and
in particular $\mathrm{\mathrm{lim}}_{1}\boldsymbol{A}_{\alpha }=0$.
\end{lemma}

\begin{proof}
For $n\in \omega $, let $B^{\left( n\right) }$ to be the set of $a\in
A^{\left( n\right) }$ such that there exists a sequence $\left( a_{i}\right)
_{i\geq n}$ with $a_{n}=a$ and, for all $k>n$, $a_{k}\in A^{\left( k\right) }
$ and $p^{\left( k-1,k\right) }\left( a_{k}\right) =a_{k-1}$. We claim that
there exists $\alpha <\omega _{1}$ such that $A_{\alpha }^{\left( n\right)
}=B^{\left( n\right) }$ for every $n\in \omega $. Let $T$ be the standard
Borel space of finite sequences $\left( x_{0},\ldots ,x_{n}\right) $ with $%
n\in \omega $ and $x_{i}\in A^{\left( i\right) }\setminus B^{\left( i\right)
}$ for $i\in \left\{ 0,1,\ldots ,n\right\} $. Let $\prec $ be the relation
on $T$ defined by setting%
\begin{equation*}
\left( x_{0},\ldots ,x_{n}\right) \prec \left( y_{0},\ldots ,y_{m}\right)
\Leftrightarrow m<n\text{, and }\forall i\leq m\text{, }x_{i}=y_{i}\text{,
and }\forall i<n\text{, }x_{i}=p^{\left( i,i+1\right) }\left( x_{i+1}\right) 
\text{.}
\end{equation*}%
We have that $\prec $ is well-founded in the sense of \cite[Appendix B]%
{kechris_classical_1995}. We also have that $\prec $ is analytic, and hence
by \cite[Theorem 31.1]{kechris_classical_1995} its rank $\rho \left( \prec
\right) $ is a countable ordinal. For $\Xi \subseteq T$ define 
\begin{equation*}
\Xi ^{\prime }=\left\{ s\in T:\exists t\in \Xi \text{, }t\prec s\right\} 
\text{.}
\end{equation*}%
Set then recursively $T_{0}=T$, $T_{\alpha +1}=T_{\alpha }^{\prime }$, and $%
T_{\lambda }=\bigcap_{\alpha <\lambda }T_{\alpha }$ for $\lambda $ limit.
Then we have that, for every ordinal $\alpha $ and $s=\left( g_{0},\ldots
,g_{m}\right) \in T$:

\begin{itemize}
\item $\rho _{\prec }\left( s\right) \geq \alpha $ if and only if $s\in
T_{\alpha }$;

\item $s\in T_{\omega \alpha }$ if and only if $g_{m}\in A_{\alpha }^{\left(
m\right) }$.
\end{itemize}

Thus, if $\lambda $ is a countable ordinal such that $\rho \left( \prec
\right) \leq \omega \lambda $, then we have that $\rho _{\prec }\left(
s\right) <\omega \lambda $ for every $s\in T$, and hence $T_{\omega \lambda
}=\varnothing $ and $A_{\lambda }^{\left( n\right) }=B^{\left( n\right) }$
for every $n\in \omega $.
\end{proof}

The same proof as \cite[Lemma 7.10]{casarosa_phantom_2025}, which applies
equally well to Polish groups that are not necessarily abelian, shows the
following:

\begin{lemma}
\label{Lemma:derived}Let $\boldsymbol{A}$ be a tower of Polish groups. For $%
\ell \leq i<\omega $ and $\alpha <\omega _{1}$ define%
\begin{equation*}
A_{\alpha }^{i}[\ell ]:=\mathrm{\mathrm{Ker}}(A_{\alpha }^{\left( i\right)
}\rightarrow A_{\alpha }^{\left( \ell \right) })\text{.}
\end{equation*}%
For fixed $\ell <\omega $ and $\alpha <\omega _{1}$, let $\boldsymbol{A}%
_{\alpha }[\ell ]$ be the tower $(A_{\alpha }^{\left( i\right) }[\ell
])_{i\geq \ell }$. Then for every $\alpha <\omega _{1}$:%
\begin{equation*}
\mathrm{lim}^{1}\boldsymbol{A}_{\alpha }=\bigcap_{\beta <\alpha
}\bigcap_{\ell <\omega }\mathrm{Ran}\left( \mathrm{lim}^{1}\boldsymbol{A}%
_{\beta }[\ell ]\rightarrow \mathrm{lim}^{1}\boldsymbol{A}\right) \text{.}
\end{equation*}%
Furthermore, if $\boldsymbol{A}$ is a tower of abelian Polish groups, then%
\begin{equation*}
\mathrm{lim}^{1}\boldsymbol{A}_{\alpha }=\mathrm{Ph}^{\alpha }\left( \mathrm{%
lim}^{1}\boldsymbol{A}\right) \text{.}
\end{equation*}
\end{lemma}

\begin{proof}
The first assertion is proved as \cite[Lemma 7.10]{casarosa_phantom_2025},
while the second assertion is the content of \cite[Theorem 7.6]%
{casarosa_phantom_2025}.
\end{proof}

\begin{definition}
Let $\boldsymbol{A}$ be a tower of \emph{countable} groups. Then on defines
the phantom length of $\boldsymbol{A}$ to be:

\begin{itemize}
\item at most $\alpha $ if and only if $\boldsymbol{A}_{\alpha }$ satisfies
the Mittag--Leffler property;

\item at most $\alpha +1/2$ if and only if $\boldsymbol{A}_{\alpha }$ is
essentially monomorphic.
\end{itemize}

Then phantom length of $\boldsymbol{A}$ is the least $\alpha \in \omega
_{1}[1/2]$ such that $\boldsymbol{A}$ has phantom length at most $\alpha $.
\end{definition}

\begin{remark}
When $\boldsymbol{A}$ is (isomorphic to) a tower of \emph{abelian }groups,
by \cite[Lemma 7.5 and Theorem 7.6]{casarosa_phantom_2025} the phantom
length of $\boldsymbol{A}$ coincides with the phantom length as in
Definition \ref{Definition:phantom-length} of the group with a Polish cover $%
\mathrm{lim}^{1}\boldsymbol{A}$.
\end{remark}

\subsection{A bound on phantom order}

Recall the notion of cofiltration of a locally compact Polish space $X$; see
also \cite[Definition 7.1]{bergfalk_definable_2024-1}. This is a sequece $%
\left( K_{n}\left( X\right) \right) _{n\in \omega }$ Notice that, if $\left(
K_{n}\left( X\right) \right) _{n\in \omega }$ of compact subsets of $X$ such
that:

\begin{itemize}
\item $K_{0}(X)=\varnothing $;

\item for every $n\in \mathbb{N}$, $K_{n}\left( X\right) $ is contained in
the interior of $K_{n+1}\left( X\right) $;

\item the sequence of interiors of $K_{n}\left( X\right) $ for $n\in \mathbb{%
N}$ form an open cover of $X$.
\end{itemize}

Notice that if $\left( K_{n}\left( X\right) \right) _{n\in \omega }$ and $%
\left( K_{n}^{\prime }\left( X\right) \right) _{n\in \omega }$ are
cofiltrations of $X$, then there exist increasing sequences $\left(
s_{n}\right) $ and $\left( t_{n}\right) $ in $\omega $ such that $%
K_{s_{n}}\left( X\right) \subseteq K_{t_{n}}\left( X\right) \subseteq
K_{s_{n+1}}\left( X\right) $ for every $n\in \omega $. In particular, $%
\left( K_{n}\left( X\right) \right) _{n\in \omega }$ and $\left(
K_{n}^{\prime }\left( X\right) \right) _{n\in \omega }$ are isomorphic in
the category of inductive sequences of compact Polish spaces.

Recall from \cite[Section 7.3]{bergfalk_definable_2024-1} the \emph{%
generalized reduced suspension} $\overline{\Sigma }\left( C\right) $ of a
compact space $C$, which is a pointed compact space. Recall also that if $C$
is a compact Polish space, and $P$ is a homotopy polyhedron, then $[C,P]$ is
countable \cite[Corollary 4.5]{bergfalk_definable_2024-1}.

If $X$ is a locally compact Polish space and $P$ is a pointed polyhedron,
then 
\begin{equation*}
\boldsymbol{G}\left( X,P\right) :=\left( [\overline{\Sigma }\left(
K_{n}(X)\right) ,P]\right) _{n\in \omega }
\end{equation*}%
is a tower of countable groups, which is independent up to natural
isomorphism from the choice of the cofiltration of $X$. It is proved in \cite%
[Theorem 7.10]{bergfalk_definable_2024} that $\mathrm{\mathrm{lim}}^{1}%
\boldsymbol{G}\left( X,P\right) $ is definably isomorphic to $\mathrm{Ph}%
^{0}[X,P]$, and the latter is the closure of $\left\{ \ast \right\} $ in $%
[X,P]$. 

\begin{lemma}
\label{Lemma:phantom-lim1}Let $X$ be a locally compact Polish space, and $P$
be a pointed polyhedron. Define $\boldsymbol{G}\left( X,P\right) $ as above.
For every countable ordinal $\alpha $, under the isomorphism%
\begin{equation*}
\mathrm{Ph}^{0}[X,P]\cong \mathrm{lim}^{1}\boldsymbol{G}\left( X,P\right) 
\text{,}
\end{equation*}%
$\mathrm{Ph}^{\alpha }[X,P]$ corresponds to $\mathrm{lim}^{1}\boldsymbol{G}%
\left( X,P\right) _{\alpha }$.
\end{lemma}

\begin{proof}
This follows from Lemma \ref{Lemma:derived} by naturality, by induction on $%
\alpha $; see also \cite[Theorem 5]{ha_higher_2003}.
\end{proof}

\begin{corollary}
Let $X$ be a locally compact Polish space, and $P$ be a pointed polyhedron.
Then there exists $\alpha <\omega _{1}$ (depending on $X$ and $P$) such that
every map $X\rightarrow P$ that is phantom of order at least $\alpha $ is
nullhomotopic.
\end{corollary}

\begin{proof}
This follows from Lemma \ref{Lemma:countable-rank} and Lemma \ref%
{Lemma:phantom-lim1}.
\end{proof}

Recall the notion of phantom length for a homogeneous space with a Polish
cover from Definition \ref{Definition:phantom-length-bound}. In view of
Lemma \ref{Lemma:phantom-lim1}, it is meaningful to consider the following:

\begin{definition}
\label{Definition:phantom-length}Let $X$ be a locally compact Polish space,
and $P$ be a pointed polyhedron. The phantom length of $[X,P]$ is the
phantom length of $\boldsymbol{G}\left( X,P\right) $.
\end{definition}

\section{Milnor\label{Section:Milnor}}

In this section we consider a general versions of the Milnor Exact Sequence
and Milnor Isomorphism, originally considered by Milnor in the context of
simplicial cohomology \cite{milnor_axiomatic_1962}; see also \cite[%
Application 3.5.9]{weibel_introduction_1995} for a version in the context of
algebraic complexes, \cite{ferry_remarks_1995} and \cite[Chapter 17]%
{mardesic_strong_2000} for its manifestation in the context of Steenrod
homology. We also obtain higher order analogues of the Milnor Isomorphism
for arbitrary phantom subgroups.

\subsection{Homotopy limits of strict complexes}

Let $A$ be a complex in $\mathbf{PAb}$. By definition, $A$ is \emph{strict}
if its differentials are strict morphisms in $\mathbf{PAb}$, i.e.,
continuous group homomorphisms \emph{with closed image}. This is equivalent
to the assertion that $\mathrm{H}_{\bullet }\left( X\right) $ is a graded
Polish group.

Suppose that $\boldsymbol{A}$ is a tower of strict chain complexes. Define%
\begin{equation*}
\mathsf{B}_{n}\left( \boldsymbol{A}\right) =\mathsf{B}_{n}\left( \mathrm{ho%
\mathrm{lim}}\boldsymbol{A}\right)
\end{equation*}%
and%
\begin{equation*}
\mathsf{Z}_{n}\left( \boldsymbol{A}\right) =\mathsf{Z}_{n}\left( \mathrm{ho%
\mathrm{lim}}\boldsymbol{A}\right)
\end{equation*}%
Notice that, for every $n\in \mathbb{Z}$, we have that%
\begin{equation*}
\left( \mathrm{ho\mathrm{lim}}\boldsymbol{A}\right) _{n}=\prod_{m\in \omega
}A_{n}^{\left( m\right) }\oplus \prod_{m_{0}<m_{1}}A_{n+1}^{\left(
m_{0}\right) }\text{.}
\end{equation*}%
The first-coordinate projection defines a continuous chain map

\begin{equation*}
\mathrm{ho\mathrm{lim}}\boldsymbol{A}\rightarrow \prod_{m\in \omega
}A^{\left( m\right) }\text{,}
\end{equation*}%
which in turn induces a continuous homomorphism%
\begin{equation*}
\mathrm{Z}_{n}\left( \boldsymbol{A}\right) \rightarrow \mathrm{Z}_{n}\left(
\prod_{m\in \omega }A^{\left( m\right) }\right) \rightarrow \prod_{m\in
\omega }\mathrm{H}_{n}(A^{\left( m\right) })\text{.}
\end{equation*}%
One can see directly that:

\begin{itemize}
\item its kernel is the closure $\overline{\mathrm{B}_{n}\left( \boldsymbol{A%
}\right) }$ of $\mathrm{B}_{n}\left( \boldsymbol{A}\right) $ inside $\mathsf{%
Z}_{n}\left( \boldsymbol{A}\right) $;

\item its image is $\mathrm{lim}_{m}\mathrm{H}_{n}\left( A^{\left( m\right)
}\right) $.
\end{itemize}

Thus, we obtain the Milnor short exact sequence%
\begin{equation}
0\rightarrow \mathrm{Ph}^{0}\mathrm{H}_{\bullet }\left( \mathrm{ho\mathrm{lim%
}}\boldsymbol{A}\right) \rightarrow \mathrm{H}_{\bullet }\left( \mathrm{ho%
\mathrm{lim}}\boldsymbol{A}\right) \rightarrow \mathrm{lim}_{m}\mathrm{H}%
_{\bullet }(A^{\left( m\right) })\rightarrow 0\text{\label{Equation:Milnor}}
\end{equation}%
which is natural in $\boldsymbol{A}$.

\subsection{Milnor isomorphism}

The Milnorm isomorphism describes $\mathrm{Ph}^{0}\mathrm{H}_{\bullet
}\left( \mathrm{ho\mathrm{lim}}\boldsymbol{A}\right) $ in terms of the
homology groups $\mathrm{H}_{\bullet }\left( \boldsymbol{A}\right) $.

\begin{lemma}
\label{Lemma:Milnor}Let $\boldsymbol{A}$ be a tower of strict complexes in $%
\mathbf{PAb}$. Then%
\begin{equation*}
\mathrm{Ph}^{0}\mathrm{H}_{\bullet }\left( \mathrm{ho\mathrm{lim}}%
\boldsymbol{A}\right) \cong \mathrm{lim}^{1}\mathrm{H}_{\bullet
+1}(A^{\left( m\right) })
\end{equation*}%
in $\mathrm{LH}\left( \mathbf{PAb}\right) $.
\end{lemma}

\begin{proof}
Suppose that $x\in \overline{\mathsf{B}_{n}\left( \boldsymbol{A}\right) }$.
Then, by definition, we have that $x_{m}\in \mathsf{B}_{n}\left( A^{\left(
m\right) }\right) $, 
\begin{equation*}
x_{m_{0},m_{1}}\in A_{n+1}^{\left( m\right) }
\end{equation*}%
, and 
\begin{equation*}
\partial x_{m_{0},m_{1}}+\left( -1\right)
^{n}(p^{(m_{0},m_{1})}(x_{m_{1}})-x_{m_{0}})=0
\end{equation*}%
for every $m_{0}<m_{1}$. Since $x_{m}\in \mathsf{B}_{n}\left( A^{\left(
m\right) }\right) $ we can choose (in a Borel fashion) $\hat{x}_{m}\in
A_{n+1}^{\left( m\right) }$ such that $\partial \hat{x}_{m}=x_{m}$. It
follows that 
\begin{equation*}
x_{m_{0},m_{1}}+\left( -1\right) ^{n}(p^{\left( m_{0},m_{1}\right) }\left( 
\hat{x}_{m_{1}}\right) -\hat{x}_{m_{0}})\in \mathsf{Z}_{n+1}(A^{\left(
m\right) })\text{.}
\end{equation*}%
This defines an element%
\begin{equation*}
x_{m_{0},m_{1}}+\left( -1\right) ^{n}(p^{\left( m_{0},m_{1}\right) }\left( 
\hat{x}_{m_{1}}\right) -\hat{x}_{m_{1}})+\mathsf{B}_{n+1}(A^{\left(
m_{0}\right) })\in \mathrm{H}_{n+1}(A^{\left( m_{0}\right) })\text{.}
\end{equation*}%
In turn, this gives an element%
\begin{equation*}
\Psi \left( x\right) :=\left( x_{m_{0},m_{1}}+\left( -1\right)
^{n}(p^{\left( m_{0},m_{1}\right) }\left( \hat{x}_{m_{1}}\right) -\hat{x}%
_{m_{0}})+\mathsf{B}_{n+1}(A^{\left( m_{0}\right) })\right)
_{m_{0}<m_{1}}\in \prod_{m_{0}<m_{1}}\mathrm{H}_{n+1}(A^{\left( m_{0}\right)
})\text{.}
\end{equation*}%
It is easily seen directly that the Borel function 
\begin{equation*}
\Psi :\overline{\mathsf{B}_{n}\left( \boldsymbol{A}\right) }\rightarrow
\prod_{m_{0}<m_{1}}\mathrm{H}_{n+1}(A^{\left( m_{0}\right) })
\end{equation*}%
induces an isomorphism%
\begin{equation*}
\mathrm{Ph}^{0}\mathrm{H}_{n}\left( \boldsymbol{A}\right) \cong \mathrm{lim}%
_{m}^{1}\mathrm{H}_{n+1}(A^{\left( m\right) })\text{.}
\end{equation*}
\end{proof}

The Milnor exact sequence together with the Milnor isomorphism produce an
exact sequence that is essentially the particular case of \cite[Theorem 5]%
{mardesic_strong_1988} in the case when $s=1$ in view of \cite[Lemma 19.34]%
{mardesic_strong_2000}; see also \cite[Theorem 17.11]{mardesic_strong_2000}.
Notice that in \cite[Theorem 5]{mardesic_strong_1988} only consider what we
call \emph{strict morphisms }between towers of complexes, while we consider
a more generous notion of morphism.

\subsection{Epimorphic towers}

Let us say that a tower $\boldsymbol{A}=\left( A^{\left( m\right) }\right) $
of strict Polish chain complexes is epimorphic if, for every $m\in \mathbb{N}
$ and $d\in \mathbb{Z}$, the component in degree $d$ of the morphisms of
complexes $A^{\left( m+1\right) }\rightarrow A^{\left( m\right) }$ is a 
\emph{surjective }continuous group homomorphism. If $\boldsymbol{A}$ is an
epimorphic tower of complexes, then we regard $\mathrm{lim}\boldsymbol{A}$
as a subcomplex of $\mathrm{holim}\boldsymbol{A}$, by identifying $a\in
\left( \mathrm{lim}\boldsymbol{A}\right) _{n}$ by the element $\hat{a}\in
\left( \mathrm{lim}\boldsymbol{A}\right) _{n}$ defined by $\hat{a}_{k}=a_{k}$
and $\hat{a}_{k,k+1}=0$ for $k\in \omega $.

For such a tower $\boldsymbol{A}$ one can obtain the Milnor exact sequence
and the Milnor isomorphism even replacing $\mathrm{ho\mathrm{lim}}%
\boldsymbol{A}$ with $\mathrm{lim}\boldsymbol{A}$, as shown in \cite[Theorem
3.5.8]{weibel_introduction_1995}, \cite[Lemma 1]{mardesic_strong_1988}, and 
\cite[Lemma 17.12]{mardesic_strong_2000}; see also \cite%
{milnor_axiomatic_1962,milnor_steenrod_1995}. Such an exact sequence,
together with the Five Lemma \cite[Section VIII.4]{mac_lane_categories_1998}%
, shows that the inclusion $\mathrm{\mathrm{lim}}\boldsymbol{A}\rightarrow 
\mathrm{ho\mathrm{lim}}\boldsymbol{A}$ induces an isomorphism in homology,
i.e., is a quasi-isomorphism.

\subsection{Recursion}

Let $\boldsymbol{A}$ be an epimorphic tower of strict complexes in $\mathbf{%
PAb}$. For $\ell \geq 0$ let $\boldsymbol{A}[\ell ]$ be the tower of strict
complexes%
\begin{equation*}
(A^{\left( m\right) }[\ell ])_{m\geq \ell }=(\mathrm{\mathrm{Ker}}(A^{\left(
m\right) }\rightarrow A^{\left( \ell \right) }))_{m\geq \ell }
\end{equation*}%
where the kernel of a morphisms of complexes is computed degree-wise. Fix $%
\ell \geq 0$. Then the inclusion%
\begin{equation*}
\mathrm{\mathrm{Ker}}(A^{\left( m\right) }\rightarrow A^{\left( \ell \right)
})\subseteq A^{\left( m\right) }
\end{equation*}%
for $m\geq \ell $ induces a morphism towers of complexes%
\begin{equation*}
\boldsymbol{A}[\ell ]\rightarrow \boldsymbol{A}\text{.}
\end{equation*}%
In turn, this induces a morphism of complexes%
\begin{equation*}
\mathrm{ho\mathrm{lim}}\boldsymbol{A}[\ell ]\rightarrow \mathrm{ho\mathrm{lim%
}}\boldsymbol{A}\text{.}
\end{equation*}%
The Milnor exact sequence for $\mathrm{H}_{\bullet }\left( \boldsymbol{A}%
\right) $ yields%
\begin{equation*}
\mathrm{Ph}^{0}\mathrm{H}_{\bullet }\left( A\right) =\bigcap_{\ell \in
\omega }\mathrm{Ran}\left( \mathrm{H}_{\bullet }\left( \boldsymbol{A}[\ell
]\right) \rightarrow \mathrm{H}_{\bullet }\left( \boldsymbol{A}\right)
\right) \text{.}
\end{equation*}%
Fix $\ell \in \omega $. Since $\boldsymbol{A}$ is an epimorphic tower, 
\begin{equation*}
\mathrm{Ran}\left( \mathrm{lim}_{m}^{1}(\mathrm{H}_{\bullet +1}(A^{\left(
m\right) })[\ell ]))\rightarrow \mathrm{lim}^{1}\mathrm{H}_{\bullet
+1}(A^{\left( \ell \right) })\right) 
\end{equation*}%
is equal to 
\begin{equation*}
\mathrm{Ran}\left( \mathrm{lim}_{m}^{1}(\mathrm{H}_{\bullet +1}(A^{\left(
m\right) }[\ell ])))\rightarrow \mathrm{lim}^{1}\mathrm{H}_{\bullet
+1}(A^{\left( \ell \right) })\right) \text{.}
\end{equation*}%
The latter corresponds to%
\begin{equation*}
\mathrm{Ran}\left( \mathrm{Ph}^{0}\mathrm{H}_{\bullet +1}(A^{\left( m\right)
}[\ell ])\rightarrow \mathrm{Ph}^{0}\mathrm{H}_{\bullet +1}(A^{\left(
m\right) })\right) 
\end{equation*}%
under the isomorphisms%
\begin{equation*}
\mathrm{Ph}^{0}\mathrm{H}_{\bullet }\left( \mathrm{\mathrm{lim}}\boldsymbol{A%
}\right) \cong \mathrm{\mathrm{lim}}_{m}^{1}\mathrm{H}_{\bullet
+1}(A^{\left( m\right) })
\end{equation*}%
and%
\begin{equation*}
\mathrm{Ph}^{0}\mathrm{H}_{\bullet }(\mathrm{\mathrm{lim}}\boldsymbol{A}%
[\ell ])\cong \mathrm{lim}_{m}^{1}\mathrm{H}_{\bullet +1}(A^{\left( m\right)
}[\ell ])\text{.}
\end{equation*}%
By the above observations, recalling the notation from Section \ref%
{Subsection:notation}, one recursively describe the phantom subgroups of $%
\mathrm{H}_{\bullet }\left( \mathrm{\mathrm{lim}}\boldsymbol{A}[\ell
]\right) $ as follows.

\begin{theorem}
\label{Theorem:phantom}Let $\boldsymbol{A}$ be an epimorphic tower of strict
complexes in\textbf{\ }$\mathbf{PAb}$. Then for every $\alpha <\omega _{1}$,%
\begin{equation*}
\mathrm{Ph}^{\alpha }\mathrm{H}_{\bullet }\left( \mathrm{lim}\boldsymbol{A}%
\right) =\bigcap_{n\in \omega }\bigcap_{\ell \in \omega }\mathrm{Ran}\left( 
\mathrm{Ph}^{\left( \alpha -1\right) _{n}}\mathrm{H}_{\bullet }\left( 
\mathrm{\mathrm{\mathrm{\mathrm{\mathrm{\mathrm{lim}}}}}}\boldsymbol{A}[\ell
]\right) \rightarrow \mathrm{H}_{\bullet }\left( \boldsymbol{A}\right)
\right) \text{.}
\end{equation*}
\end{theorem}

\section{Phantom homology and cohomology\label{Section:phantom}}

In this section we apply the Universal Coefficient Theorem and the Milnor
Exact Sequence and Milnor Isomorphisms to describe the phantom subgroups of
Steenrod homology and \v{C}ech cohomology, and relate them with phantom maps.

\subsection{Steenrod property}

An abelian Polish group $G$ is Steenrod or has the \emph{Steenrod property }%
(also called \emph{division-closure property} by Steenrod \cite%
{steenrod_universal_1936} ad Lefshetz \cite{lefschetz_algebraic_1942}) if,
for every $k\in \mathbb{N}$,%
\begin{equation*}
kG=\left\{ kg:g\in G\right\}
\end{equation*}%
is a closed subgroup of $G$.

It is clear that if $G$ is compact, then $G$ has the Steenrod property.\ The
same holds for vector groups, and discrete groups, and hence for all
products of abelian Lie groups. The following example of a locally compact
Polish abelian group that does \emph{not }have the Steenrod property was
pointed out by Braunling. Let $G:=\mathbb{Q}_{p}$ be the additive group of
the ring of $p$-adic numbers and $H=\mathbb{Z}_{p}\subseteq G$ be the ring
of $p$-adic integers. Define $L$ to be the restricted product of copies of $G
$ with respect to the open subgroup $H$. By definition, $L$ is the subgroup
of $G^{\omega }$ comprising the sequences $\left( x_{n}\right) $ of elements
of $G$ that belong eventually to $H$.\ The topology on $L$ is obtained by
declaring $H^{\omega }$ with the product topology to be an open subgroup, so
that%
\begin{equation*}
0\rightarrow H^{\omega }\rightarrow L\rightarrow \left( G/H\right) ^{\left(
\omega \right) }\rightarrow 0
\end{equation*}%
is a short exact sequence. Then for every $n\in \mathbb{N}$, $p^{n}L$ is a
dense proper subgroup of $L$. A related example is provided by the additive
group of the (finite) adele over $\mathbb{Q}$ \cite[Chapter IV]%
{weil_basic_1995}.

\subsection{Homology\label{Subsection:homology}}

The (Steenrod) homology $\mathrm{H}^{\bullet }\left( X;G\right) $ of a
compact Polish space $X$ with coefficients in an abelian group with a Polish
cover $G$ is obtained as in Section \ref{Section:homology-spaces} by letting 
$\mathcal{M}$ be the countably complete abelian category $\mathrm{LH}\left( 
\mathbf{PAb}\right) $ of abelian groups with a Polish cover. As a particular
instance of Theorem \ref{Theorem:UCT-homology-spaces}, we obtain the
following Universal Coefficient Theorem for homology of compact Polish
spaces.

\begin{theorem}
\label{Theorem:UCT-homology}Suppose that $X$ is a compact metrizable space,
and let $G$ be a group with a Polish cover. Then we have a natural short
exact sequence of groups with a Polish cover%
\begin{equation*}
0\rightarrow \mathrm{Ext}\left( \mathrm{H}^{\bullet +1}\left( X\right)
,G\right) \rightarrow \mathrm{H}_{\bullet }\left( X;G\right) \rightarrow 
\mathrm{Hom}\left( \mathrm{H}^{\bullet }\left( X\right) ,G\right)
\rightarrow 0
\end{equation*}%
which splits (unnaturally).
\end{theorem}

It follows from the structure of group homomorphisms between free finitely
generated abelian groups given by the Smit\`{h} Normal Form for integer
matrices \cite{newman_smith_1997} that, for every complex of free
finitely-generated abelian groups $A$, if $G$ is an abelian Polish group
with the Steenrod property, then $\mathrm{Hom}\left( A,G\right) $ is a \emph{%
strict} complex in $\mathbf{K}\left( \mathbf{PAb}\right) $. (In other words,
the differential maps in $\mathrm{Hom}\left( A,G\right) $ have closed
image.) It follows that, for a Polish group $G$ with the Steenrod property
and a countable abelian group $A$, \textrm{Ph}$^{0}\mathrm{Ext}\left(
A,G\right) $ is the closed subgroup $\mathrm{PExt}\left( A,G\right) $
corresponding to \emph{pure }extensions. By definition, this is the kernel
of the canonical homomorphism%
\begin{equation*}
\mathrm{Ext}\left( A,G\right) \rightarrow \mathrm{lim}_{B}\mathrm{Ext}\left(
B,G\right)
\end{equation*}%
where $B$ ranges among the \emph{finitely-generated} subgroups of $A$
ordered by inclusion.

Suppose that $X$ is a compact Polish space. Then one can associate with $X$
a tower $\left( X^{\left( m\right) }\right) $ of compact polyhedra, called 
\emph{polyhedral resolution} of $X$, such that $X$ is homeomorphic to the
inverse limit $\mathrm{lim}_{m}X^{\left( m\right) }$; see \cite[Theorem 7,
page 61 and Corollary 6, page 67]{mardesic_shape_1982}. One can define such
a tower in terms of a cofinal sequence $\mathcal{U}^{X}=\left( \mathcal{U}%
_{m}^{X}\right) _{m\in \omega }$ of finite open covers for $X$, setting $%
X^{\left( m\right) }:=\left\vert N\left( \mathcal{U}_{m}^{X}\right)
\right\vert $ for $m\in \omega $; see \cite[Section IX.9]%
{eilenberg_foundations_1952}. Such a polyhedral resolution yields a
surjective Borel-definable group homomorphism%
\begin{equation*}
\mathrm{H}_{\bullet }\left( X;G\right) \rightarrow \mathrm{\mathrm{lim}}_{m}%
\mathrm{H}_{\bullet }\left( X_{m};G\right) 
\end{equation*}%
where $\mathrm{\mathrm{lim}}_{m}\ \mathrm{H}_{\bullet }\left( X_{m};G\right) 
$ is a Polish group. Using this fact, the Universal Coefficient Sequence,
and its naturality in $X$, we obtain the following.

\begin{corollary}
\label{Corollary:UCT-homology}Let $X$ be a compact Polish space, and let $G$
be an abelian Polish group with the Steenrod property. Suppose that $\left(
X^{\left( m\right) }\right) $ is a polyhedral resolution of $X$. Then%
\begin{equation*}
\mathrm{Ph}^{0}\mathrm{H}_{\bullet }\left( X;G\right) =\mathrm{Ker}\left( 
\mathrm{H}_{\bullet }\left( X;G\right) \rightarrow \mathrm{\mathrm{lim}}_{m}%
\mathrm{H}_{\bullet }\left( X_{m};G\right) \right) \cong \mathrm{PExt}\left( 
\mathrm{H}^{\bullet +1}\left( X\right) ,G\right) \text{.}
\end{equation*}
\end{corollary}

In the context of Corollary \ref{Corollary:UCT-homology}, we define%
\begin{equation*}
\mathrm{H}_{\bullet }^{\mathrm{w}}\left( X;G\right) :=\mathrm{\mathrm{lim}}%
_{m}\mathrm{H}_{\bullet }\left( X_{m};G\right) \cong \frac{\mathrm{H}%
_{\bullet }\left( X;G\right) }{\mathrm{Ph}^{0}\mathrm{H}_{\bullet }\left(
X;G\right) }\text{.}
\end{equation*}%
The Milnor isomorphism in this case yields%
\begin{equation*}
\mathrm{Ph}^{0}\mathrm{H}_{\bullet }\left( X;G\right) \cong \mathrm{lim}%
_{m}^{1}\ \mathrm{H}_{\bullet +1}\left( X_{m};G\right) \text{.}
\end{equation*}

\subsection{Cohomology}

Analogous considerations as those from Section \ref{Subsection:homology}
pertaining to homology of compact Polish spaces apply to cohomology of
homotopy polyhedra. As a particular instance of Theorem \ref%
{Theorem:UCT-cohomology-polyhedra} one obtains the following Universal
Coefficient Theorem for cohomology:

\begin{theorem}
\label{Theorem:UCT-cohomology}Suppose that $Y$ is a homotopy polyhedron, and
let $G$ be an abelian group with a Polish cover. Then we have a natural
short exact sequence of groups with a Polish cover%
\begin{equation*}
0\rightarrow \mathrm{Ext}\left( \mathrm{H}_{\bullet -1}\left( Y\right)
,G\right) \rightarrow \mathrm{H}^{\bullet }\left( Y;G\right) \rightarrow 
\mathrm{Hom}\left( \mathrm{H}_{\bullet }\left( Y\right) ,G\right)
\rightarrow 0
\end{equation*}%
which splits (unnaturally).
\end{theorem}

If $Y$ is a homotopy polyhedron, then up to homotopy we can write $Y$ as $%
\mathrm{hoco\mathrm{lim}}_{m}Y_{m}$ where $Y_{m}$ are homotopy equivalent to 
\emph{finite} polyhedra. By homotopy-invariance, we have a canonical
surjective Borel-definable homomorphism%
\begin{equation*}
\mathrm{H}^{\bullet }\left( Y;G\right) \rightarrow \mathrm{\mathrm{lim}}%
_{m}\ \mathrm{H}^{\bullet }\left( Y_{m};G\right) \text{.}
\end{equation*}%
As in the case of homology for compact spaces, we obtain from the Universal
Coefficient Theorem the following consequence.

\begin{corollary}
\label{Corollary:UCT-cohomology}Suppose that $Y$ is a homotopy polyhedron,
and let $G$ be an abelian Polish group with the Steenrod property. Suppose
that $Y$ is homotopy equivalent to $\mathrm{hocolim}_{m}Y_{m}$, where $Y_{m}$
is homotopy equivalent to a finite polyhedron. Then%
\begin{equation*}
\mathrm{Ph}^{0}\mathrm{H}^{\bullet }\left( Y;G\right) =\mathrm{Ker}\left( 
\mathrm{H}^{\bullet }\left( Y;G\right) \rightarrow \mathrm{\mathrm{lim}}_{m}%
\mathrm{H}^{\bullet }\left( Y_{m};G\right) \right) \cong \mathrm{PExt}\left(
H_{\bullet -1}\left( Y\right) ,G\right) \text{.}
\end{equation*}
\end{corollary}

In the context of Corollary \ref{Corollary:UCT-cohomology}, we define%
\begin{equation*}
\mathrm{H}_{\mathrm{w}}^{\bullet }\left( Y;G\right) =\mathrm{\mathrm{lim}}%
_{m}\mathrm{H}^{\bullet }\left( Y_{m};G\right) \cong \frac{\mathrm{H}%
^{\bullet }\left( Y;G\right) }{\mathrm{Ph}^{0}\mathrm{H}^{\bullet }\left(
Y;G\right) }\text{.}
\end{equation*}%
The Milnor isomorphism in this case yields%
\begin{equation*}
\mathrm{Ph}^{0}\mathrm{H}^{\bullet }\left( Y;G\right) \cong \mathrm{lim}%
_{m}^{1}\ \mathrm{H}^{\bullet -1}\left( Y_{m};G\right) \text{.}
\end{equation*}

By way of Theorem \ref{Theorem:phantom}, one can generalize Corollary \ref%
{Corollary:UCT-cohomology} to higher order phantom subgroups, referring
again to the notation from Section \ref{Subsection:notation}:

\begin{theorem}
Let $Y$ be a homotopy polyhedron, and $G$ be a Polish abelian group with the
Steenrod property. Suppose that $\left( Y_{m}\right) $ is a polyhedral
approximation of $Y$. Then for every $\alpha <\omega _{1}$ and $d\in \mathbb{%
N}$:%
\begin{eqnarray*}
\mathrm{Ph}^{\alpha }\mathrm{H}^{d}\left( Y;G\right) &=&\bigcap_{n<\omega
}\bigcap_{m<\omega }\mathrm{Ran}\left( \mathrm{Ph}^{\left( \alpha -1\right)
_{n}}\mathrm{H}^{d}\left( Y,Y_{m};G\right) \rightarrow \mathrm{H}^{d}\left(
Y;G\right) \right) \\
&=&\bigcap_{n<\omega }\bigcap_{m<\omega }\mathrm{Ran}\left( \mathrm{Ph}%
^{\left( \alpha -1\right) _{n}}\mathrm{H}^{d}\left( Y/Y_{m};G\right)
\rightarrow \mathrm{H}^{d}\left( Y;G\right) \right)
\end{eqnarray*}%
where $\mathrm{H}^{d}\left( Y,Y_{m};G\right) $ is the cohomology of the pair 
$\left( Y,Y_{m}\right) $.
\end{theorem}

Theorem \ref{Theorem:phantom} also follows from \cite[Theorem 8.3]%
{casarosa_phantom_2025} applied to the cohomological functor $\mathrm{H}%
^{0}\left( -\right) $ on the derived category of countable abelian groups.

As a particular case of Theorem \ref{Theorem:UCT-cohomology}, one obtains
that when $Y$ is a homotopy polyhedron such that $\mathrm{H}_{d}\left(
Y\right) $ and $\mathrm{H}_{d-1}\left( Y\right) $ are finitely-generated,
then so is $\mathrm{H}^{d}\left( Y\right) $. If furthermore $\mathrm{H}%
_{d}\left( Y\right) $ and $\mathrm{H}_{d-1}\left( Y\right) $ are free and
finitely-generated, then $\mathrm{H}^{d}\left( Y\right) $ is naturally
isomorphic to $\mathrm{H}_{d}\left( Y\right) ^{\bot }$. Here, we denote by $%
A\mapsto A^{\bot }$ the functor that assigns to a free finitely-generated
abelian group $A$ the group $\mathrm{Hom}\left( A,\mathbb{Z}\right) $.
Similar considerations apply to compact Polish spaces, reversing the roles
of homology and cohomology.

\subsection{Phantom extensions}

Let $C$ and $A$ be countable abelian groups. Then the phantom subgroup $%
\mathrm{Ph}^{\alpha }\mathrm{Ext}\left( C,A\right) $ is characterized in 
\cite[Section 9.2]{casarosa_phantom_2025} as the subgroup parametrizing the
extensions of $C$ by $A$ that are \emph{phantom }(or \emph{pure}) \emph{of
order }at least $\alpha $. As a consequence of the Universal Coefficient
Theorem, and the fact that the phantom subgroups are subfunctors of the
identity on $\mathrm{LH}\left( \mathbf{PAb}\right) $, we obtain the
following:

\begin{theorem}
Let $Y$ be a homotopy polyhedron and $A$ be a countable abelian group. For $%
d\geq 1$, $\mathrm{Ph}^{\alpha }\mathrm{H}^{d}\left( Y;A\right) $ is
isomorphic to the group $\mathrm{Ph}^{\alpha }\mathrm{Ext}\left( \mathrm{H}%
_{d-1}\left( Y\right) ,A\right) $ parametrizing phantom extensions of $%
\mathrm{H}_{d-1}\left( Y\right) $ by $A$ of order at least $\alpha $.
\end{theorem}

If $C$ is a countable abelian group, then a family of quotients $\partial
_{\alpha }C$ for $\alpha <\omega _{1}$ of $C$ has been defined in \cite[%
Chapter 11]{casarosa_phantom_2025} by recursion. Then one has that the
quotient map $C\rightarrow \partial _{\alpha }C$ induces an isomorphism%
\begin{equation*}
\mathrm{Ext}\left( \partial _{\alpha }C,\mathbb{Z}\right) \cong \mathrm{Ph}%
^{\alpha }\mathrm{Ext}\left( C,\mathbb{Z}\right) \text{;}
\end{equation*}%
see \cite[Theorem 11.5]{casarosa_phantom_2025}.

\begin{definition}
\label{Definition:projective-length}Let $C,A$ be countable abelian groups.
The projective $A$-length is the phantom length of $\mathrm{Ext}\left(
C,A\right) $.
\end{definition}

In particular, the projective $\mathbb{Z}$-length of $C$ is the phantom
length of $\mathrm{Ext}\left( C,\mathbb{Z}\right) $. It follows from the
above remarks that the projective $\mathbb{Z}$-length of $C$ is at most $%
\alpha $ if and only if $\partial _{\alpha }C$ is a free abelian group, and
at most $\alpha +1/2$ if and only if $\partial _{\alpha }C$ is the sum of a
free abelian group and a finite-rank torsion-free abelian group; see \cite[%
Chapter 11]{casarosa_phantom_2025}.

\subsection{Steenrod duality\label{SubsectionDuality}}

Steenrod homology was introduced by Steenrod in \cite{steenrod_regular_1940}
as the correct \emph{dual} of \v{C}ech cohomology theory. Indeed, one of the
main results of \cite{steenrod_regular_1940} is the following \emph{Steenrod
Duality Theorem}, generalizing Alexander duality \cite[Theorem 2.15]%
{prasolov_elements_2007}, and later generalized by Sitnikov \cite%
{sitnikov_duality_1951}; see also \cite[Corollary 11.21]%
{massey_homology_1978}.

\begin{theorem}[Steenrod]
\label{Theorem:Steenrod-duality}Fix $n\in \mathbb{N}$, and $X\subseteq
S^{n+1}$ closed. If $k\in \left\{ 0,1,\ldots ,n\right\} $, then $\mathrm{H}%
^{k}\left( S^{n+1}\setminus X\right) \cong \mathrm{\widetilde{H}}%
_{n-k}\left( X\right) $, where the isomorphism is natural with respect to
the inclusion maps.
\end{theorem}

Here, $\mathrm{\widetilde{H}}_{\bullet }\left( X\right) $ denotes the \emph{%
reduced }Steenrod homology, which is by definition equal to%
\begin{equation*}
\mathrm{\mathrm{Ker}}\left( \mathrm{H}_{\bullet }\left( X\right) \rightarrow
H_{\bullet }\left( \left\{ \ast \right\} \right) \right)
\end{equation*}%
where one considers the morphism 
\begin{equation*}
\mathrm{H}_{\bullet }\left( X\right) \rightarrow \mathrm{H}_{\bullet }\left(
\left\{ \ast \right\} \right)
\end{equation*}%
associated with the canonical (constant) map $X\rightarrow \left\{ \ast
\right\} $.

We prove here the natural \emph{definable version }of the Steenrod Duality
Theorem. Suppose that $K$ is a finite simplicial complex, and that $L$ is a
simplicial subcomplex of $K$. If $\left\vert K\right\vert $ is a geometric
realization of $K$, then the union of the closed subsets of $\left\vert
K\right\vert $ corresponding to the simplices of $L$ is a geometric
realization of $L$. Notice that the inclusion of $\left\vert L\right\vert $
into $\left\vert K\right\vert $ is a cofibration \cite[Proposition 3.2.4]%
{arkowitz_introduction_2011}, and hence $\left\vert L\right\vert $ is a
deformation retract of a neighborhood of $\left\vert L\right\vert $ in $%
\left\vert K\right\vert $. The \emph{closed combinatorial neighborhood }$%
U_{K}\left( L\right) $ of $L$ in $K$ is the collection of all simplices of $%
K $ together with their faces that have some simplex of $L$ as a face. The 
\emph{open combinatorial neighborhood }$U_{K}^{\circ }\left( L\right) $ of $%
\left\vert L\right\vert $ on $\left\vert K\right\vert $ is the union of the
relative interiors of the closed subsets of $\left\vert K\right\vert $
corresponding to the simplices of $U_{K}\left( L\right) $. Notice that the
complement of $U_{K}^{\circ }\left( L\right) $ in $\left\vert K\right\vert $
can be seen as the topological realization of the simplicial complex
consisting of the simplices of $K$ that do not have a simplex of $L$ as a
face.

\begin{theorem}
\label{Theorem:definable-Steenrod-duality}Fix $n\in \mathbb{N}$ and $%
X\subseteq S^{n+1}$ closed. Then:

\begin{itemize}
\item $\mathrm{H}_{\mathrm{w}}^{k}\left( Y\right) $ and $\mathrm{\widetilde{H%
}}_{n-k}^{\mathrm{w}}\left( X\right) $ are isomorphic for $k\in \left\{
1,2,\ldots ,n\right\} $;

\item $\mathrm{Ph}^{0}\mathrm{H}^{k}(S^{n+1}\setminus X)$ and $\mathrm{Ph}%
^{0}\mathrm{\widetilde{H}}_{n-k}\left( X\right) $ are isomorphic for $k\in
\left\{ 2,3,\ldots ,n\right\} $.
\end{itemize}
\end{theorem}

\begin{proof}
Fix $k\in \left\{ 0,\ldots ,n\right\} $. Set $Y:=S^{n+1}\setminus X$. The
proof is inspired by the proof of the Alexander--Pontryagin duality theorem
from \cite[Theorem 5.7]{prasolov_elements_2007}. Fix an increasing sequence $%
\left( Y_{m}\right) _{m\in \omega }$ of compact subsets of $Y$ with union
equal to $Y$. In the proof, we will identify triangulations of $S^{n+1}$
with the corresponding (abstract) simplicial complexes. We will also
identify simplices with the corresponding closed subsets of $S^{n+1}$.

We can regard $S^{n+1}$ be the boundary of the standard geometric simplex $%
\Delta ^{n+2}$. This gives a triangulation $T_{-1}$ of $S^{n+1}$. One can
define recursively triangulations $T_{m}$ of $S^{n+1}$ and simplicial
subcomplexes $L_{m}$ and $K_{m}$ of $T_{m}$ for $m\in \omega $ such that:

\begin{itemize}
\item $T_{m+1}$ is obtained from $T_{m}$ by taking barycentric subdivisions;

\item $L_{m}$ is the set of simplices in $T_{m}$ that have nonempty
intersection with $Y_{m}$;

\item the open combinatorial neighborhood $U_{T_{m}}^{\circ }(L_{m})$ is
contained in $Y$ and has $L_{m}$ as a deformation retract;

\item $K_{m}$ is the set of simplices in $T_{m}$ that do not have a simplex
of $L_{m}$ as a face.
\end{itemize}

As $K_{m}$ is the complement of $U_{T_{m}}^{\circ }(L_{m})$, we have that $%
X\subseteq K_{m}$ for every $m\in \omega $. Thus, $X$ is homeomorphic to the
inverse limit of the sequence $\left( K_{m},\iota _{m}\right) $ where $\iota
_{m}:K_{m+1}\rightarrow K_{m}$ is the inclusion map. Furthermore, we have
that $Y$ is the union of the increasing sequence of closed subsets $%
(L_{m})_{m\in \omega }$. By the Steenrod (or Alexander) duality theorem, we
have that, for $1\leq k\leq n$,%
\begin{equation*}
\mathrm{H}^{k}(L_{m})\cong \mathrm{H}^{k}\left( U_{T_{m}}^{\circ
}(L_{m})\right) \cong \mathrm{\widetilde{H}}_{n-k}(K_{m})\text{.}
\end{equation*}%
Since such isomorphisms are natural with respect to inclusion maps, they
induce an isomorphism between the towers $\left( \mathrm{H}%
^{k}(L_{m})\right) _{m\in \omega }$ and $(\mathrm{\widetilde{H}}%
_{n-k}(K_{m}))_{m\in \omega }$. Hence, we have that, for $1\leq k\leq n$,%
\begin{equation*}
\mathrm{H}_{\mathrm{w}}^{k}\left( Y\right) \cong \mathrm{lim}_{m}{}\mathrm{H}%
^{k}(L_{m})\cong \mathrm{lim}_{m}{}\mathrm{\widetilde{H}}_{n-k}(K_{m})\cong
{}\mathrm{\widetilde{H}}_{n-k}^{\mathrm{w}}\left( X\right) \text{.}
\end{equation*}%
Similarly, for $1\leq k\leq n-1$,%
\begin{equation*}
\mathrm{Ph}^{0}\mathrm{H}^{k+1}\left( Y\right) \cong \mathrm{lim}_{m}^{1}{}%
\mathrm{H}^{k}(L_{m})\cong \mathrm{lim}_{m}^{1}{}\mathrm{\widetilde{H}}%
_{n-k}(K_{m})\cong \mathrm{Ph}^{0}\mathrm{\widetilde{H}}_{n-\left(
k+1\right) }\left( X\right) \text{.}
\end{equation*}%
This concludes the proof.
\end{proof}

\subsection{Homotopical description of cohomology}

Let $G$ be a countable abelian group, and $n\geq 1$. Let $K\left( G,n\right) 
$ be an \emph{Eilenberg--MacLane space }of type $\left( G,n\right) $ \cite[%
Definition 5.4]{bergfalk_definable_2024-1}; see also \cite[Definition 2.5.7]%
{arkowitz_introduction_2011}. For a locally compact Polish space $X$, $%
[X,K\left( G,n\right) ]$ is a group, being $K\left( G,n\right) $ a
polyhedral $H$-group. \emph{Huber's Theorem }\cite{huber_homotopical_1961}
establishes a natural group isomorphism%
\begin{equation*}
\mathrm{H}^{n}\left( X;G\right) \cong \lbrack X,K\left( G,n\right) ]
\end{equation*}%
for a locally compact Polish spacer $X$. Such an isomorphism is in fact
Borel-definable, i.e., induced by a Borel function from continuous functions 
$X\rightarrow K\left( G,n\right) $ to the space of cycles $\mathrm{Z}%
^{n}\left( X;G\right) $ \cite[Theorem 5.7]{bergfalk_definable_2024-1},
whence an isomorphism of definable groups. Via this isomorphism, $\mathrm{Ph}%
^{0}\mathrm{H}^{n}\left( X;G\right) $ corresponds to $\mathrm{Ph}%
^{0}[X,K\left( G,n\right) ]$ and, more generally, $\mathrm{Ph}^{\alpha }%
\mathrm{H}^{n}\left( X;G\right) $ to $\mathrm{Ph}^{\alpha }[X,K\left(
G,n\right) ]$ for every $\alpha <\omega _{1}$.

Combining this with the Universal Coefficient Theorem, we get a
correspondence%
\begin{equation*}
\mathrm{Ph}^{\alpha }[X,K\left( G,n\right) ]\leftrightarrow \mathrm{Ph}%
^{\alpha }\mathrm{Ext}\left( \mathrm{H}_{n-1}\left( X\right) ,G\right) 
\end{equation*}%
between homotopy class of phantom maps $X\rightarrow K\left( G,n\right) $ of
order at least $\alpha $ and isomorphism classes of extensions of $\mathrm{H}%
_{n-1}\left( X\right) $ by $G$ that are phantom of order at least $\alpha $.

By Lemma \ref{Lemma:phantom-lim1}, one also has a Borel-definable bijection%
\begin{equation*}
\mathrm{Ph}^{0}[X,K\left( G,n\right) ]\leftrightarrow \mathrm{lim}^{1}{}%
\boldsymbol{A}
\end{equation*}%
where $\boldsymbol{A}$ is the tower of countable \emph{abelian }groups $%
\boldsymbol{G}\left( X,K\left( G,n\right) \right) $. It follows that the
following elements of $\omega _{1}[1/2]$ are equal:

\begin{itemize}
\item the phantom length of $[X,K\left( G,n\right) ]$ as in Definition \ref%
{Definition:phantom-length};

\item the phantom length of $\mathrm{H}^{n}\left( X;G\right) $;

\item the phantom length of $\mathrm{Ext}\left( \mathrm{H}_{n-1}\left(
X\right) ,G\right) $;

\item (when $G$ is countable) the projective $G$-length of $\mathrm{H}%
_{n-1}\left( X\right) $.
\end{itemize}

\subsection{The Hopf Classification Theorem}

Fix $n\geq 1$ and let $S^{n}$ be the $n$-dimensional sphere. Then $\mathrm{H}%
^{n}\left( S^{n}\right) \cong \mathbb{Z}$. Let $\iota :S^{n}\rightarrow
K\left( \mathbb{Z},n\right) $ be a map whose homotopy class is a generator
for $[S^{n},K\left( \mathbb{Z},n\right) ]\cong \mathrm{H}^{n}\left(
S^{n}\right) $.

Suppose that $Q$ is a homotopy polyhedron with $H^{k}\left( Q\right) =0$ for 
$k>n$. The \emph{Hopf Classification Theorem} \cite[Chapter VII]%
{hu_homotopy_1959} asserts that the function%
\begin{equation*}
\lbrack Q,S^{n}]\rightarrow \lbrack Q,K\left( \mathbb{Z},n\right) ]\cong 
\mathrm{H}^{n}\left( Q\right) 
\end{equation*}%
mapping the homotopy class of $f:Q\rightarrow S^{n}$ to the homotopy class
of $\iota \circ f:Q\rightarrow K\left( \mathbb{Z},n\right) $, is a
bijection. Again, via this bijection, $\mathrm{Ph}^{\alpha }[Q,S^{n}]$
corresponds to 
\begin{equation*}
\mathrm{Ph}^{\alpha }\mathrm{H}^{n}\left( Q\right) \cong \mathrm{Ph}^{\alpha
}\mathrm{Ext}\left( \mathrm{H}_{n-1}\left( Q\right) ,\mathbb{Z}\right) 
\end{equation*}%
for every $\alpha <\omega _{1}$.

If $\gamma _{n}\in \mathrm{H}_{n}\left( S^{n}\right) \cong \mathbb{Z}$ is a
generator, then the \emph{Hurewicz homomorphism}%
\begin{equation*}
\pi _{n}\left( Y\right) \rightarrow \mathrm{H}_{n}\left( Y\right) 
\end{equation*}%
is defined by mapping an element of $\pi _{n}\left( Y\right) =[S_{n},Y]$
represented by a map $f:S_{n}\rightarrow Y$ to $\mathrm{H}_{n}\left(
f\right) \left( \gamma _{n}\right) $ \cite[Section 2.4]%
{arkowitz_introduction_2011}. If $n\geq 1$ and $Y$ is a simply connected
space (i.e., connected with trivial fundamental group) such that%
\begin{equation*}
H_{m}\left( Y\right) =0\text{ for }1<m<n\text{,}
\end{equation*}%
then the Hurewicz homomorphism $\pi _{n}\left( Y\right) \rightarrow \mathrm{H%
}_{n}\left( Y\right) $ is an isomorphism \cite[Theorem 5.6.17]%
{arkowitz_introduction_2011}.

\section{Higher order Borsuk--Eilenberg classification\label{Section:higher}}

In this section, we consider natural higher order analogues of (concrete)
solenoids as compact subspaces of $\mathbb{R}^{3}$. We then apply the
Universal Coefficient Theorem to classify the continuous maps from their
complements to a $2$-dimensional sphere up to homotopy. This classification
problem was considered by Borsuk and Eilenberg in 1936 in the case of
standard $p$-adic solenoids \cite{borsuk_uber_1936}.\ Recently, this problem
was studied from the perspective of Borel complexity theory in \cite%
{bergfalk_definable_2024-1}. In this section we show in particular that,
when one replaces $p$-adic solenoids with higher order solenoids, the
corresponding classification problem can have arbitrarily high complexity.

\subsection{Solenoid complements}

Recall that a $d$-dimensional solenoid is an indecomposable continuum, i.e.\
a compact connected topological space, and in particular a compact Polish
space \cite{vietoris_uber_1927,van_dantzig_uber_1930}; see also \cite%
{mccord_inverse_1965,aarts_classification_1991}. Such spaces can be
described, up to homeomorphism, as Pontryagin duals of infinitely generated
torsion-free abelian groups of rank $d$ \cite%
{bognar_embedding_1988,bognar_embedding_1988-1}.

Fix an integer $p\geq 2$ and let $\Sigma _{p}$ denote the $p$-adic solenoid,
namely the Pontryagin dual of $\mathbb{Z}[1/p]$, the additive group of
rational numbers whose denominators are powers of $p$.

We recall the construction of a concrete embedding of $\Sigma _{p}$ as a
closed subspace of $\mathbb{R}^{3}$. Let $\left( W_{k}\right) $ be a
decreasing sequence of compact subsets of $\mathbb{R}^{3}$ such that, for
every $k\in \omega $:

\begin{enumerate}
\item $W_{k}$ is homeomorphic to $\mathbb{T}\times D$, where $D$ is the $2$%
-dimensional disk;

\item the inclusion $W_{k+1}\rightarrow W_{k}$ has degree $p$;
\end{enumerate}

see \cite[Exercise VIII.E]{eilenberg_foundations_1952}. Then one sets%
\begin{equation*}
\Sigma _{p}:=\bigcap_{k\in \omega }W_{k}\subseteq \mathbb{R}^{3}\text{.}
\end{equation*}%
Explicit constructions can be found in \cite%
{hubbard_henon_1994,bognar_embedding_1988,bognar_embedding_1988-1,conner_geometry_2015,jiang_tame_2011,clark_embedding_2004}%
. Concrete solenoids arising as attractors of dynamical systems are usually
called Smale solenoids or Smale--Williams solenoids \cite%
{yu_smale_2019,ma_realization_2007,ma_genus_2011,smale_differentiable_1967};
see also \cite[Chapter VI]{stewart_does_2002}.

It follows from this and the Universal Coefficients Theorem that%
\begin{equation*}
\mathrm{H}^{2}\left( S^{3}\setminus \Sigma _{p}\left( h\right) \right) \cong 
\mathrm{Ext}\left( \mathbb{Z}[1/p],\mathbb{Z}\right)
\end{equation*}%
and%
\begin{equation*}
\mathrm{H}^{d}\left( S^{3}\setminus \Sigma _{p}\left( h\right) \right) =0
\end{equation*}%
for $d>2$. The Hopf Classification Theorem thus produces a definable
bijection 
\begin{equation*}
\lbrack S^{3}\setminus \Sigma _{p}(h),S^{2}]\longleftrightarrow \mathrm{Ext}(%
\mathbb{Z}[1/p],\mathbb{Z})\text{.}
\end{equation*}

\subsection{Handlebodies}

We now extend the construction of solenoids to higher order versions. Higher
genus solenoids have also been considered in \cite{ma_genus_2011} in the
context dynamical systems. Recall that if $G$ is a polyhedron embedded in a $%
3$-manifold, a \emph{regular neighborhood} of $G$ is a closed neighborhood $%
N $ of $G$, homeomorphic to a $3$-manifold (possibly with boundary), such
that $G$ is a spine of $N$, i.e.\ $N$ deformation retracts onto $G$. A \emph{%
handlebody of genus $g$} is a $3$-manifold obtained as a regular
neighborhood of a connected simple graph of genus $g$ embedded in $\mathbb{R}%
^{3}$. A handlebody of genus $g$ can be realized as a regular neighborhood
of a wedge of $g$ circles embedded in $\mathbb{R}^{3}$. 

Fix $d\geq 1$. Let $W$ is a handlebody of genus $d$ with spine $S$. In
particular, 
\begin{equation*}
\widetilde{\mathrm{H}}_{k}(W)\cong \left\{ 
\begin{array}{ll}
\mathbb{Z}^{d} & \text{if }k=1\text{,} \\ 
0 & \text{otherwise.}%
\end{array}%
\right. \ .
\end{equation*}%
Likewise, by the Hurewicz isomorphism%
\begin{equation*}
\mathbb{Z}^{d}\cong \pi _{1}\left( W\right) \cong \mathrm{H}_{1}\left(
W\right) \text{.}
\end{equation*}%
Fix $e_{1},\ldots ,e_{d}$ be circles in $S$ that yield generators $%
[e_{1}],\ldots ,[e_{d}]$ of $\pi _{1}\left( W\right) $.

Suppose now that $k\geq d$ and $\left( a_{ij}\right) $ is a $k\times d$
integer matrix, representing a homomorphism $\varphi :\mathbb{Z}%
^{k}\rightarrow \mathbb{Z}^{d}$. As in the genus-one case, one can find a
handlebody $V\subseteq W$ of genus $k$ and spine $T$ and pairwise disjoint
circles $f_{1},\ldots ,f_{k}$ in $T$ such that:

\begin{enumerate}
\item $[f_{1}],\ldots ,[f_{k}]$ are generators of $\pi _{1}\left( V\right) $;

\item for $1\leq i\leq k$, the inclusion $f_{i}\rightarrow W$ induces the
homomorphism%
\begin{equation*}
\pi _{1}\left( f_{i}\right) \rightarrow \pi _{1}\left( W\right) \text{, }%
[f_{i}]\mapsto a_{i1}[e_{1}]+\cdots +a_{id}[e_{d}]
\end{equation*}

\item the inclusion $V\rightarrow W$ induces the homomorphism%
\begin{equation*}
\pi _{1}\left( V\right) \rightarrow \pi _{1}\left( W\right) \text{, }%
[f_{i}]\mapsto a_{i1}[e_{1}]+\cdots +a_{id}[e_{d}]\text{ for }1\leq i\leq k%
\text{.}
\end{equation*}
\end{enumerate}

It follows via the Hurewicz isomorphism that the inclusion $V\rightarrow W$
induces in homology of degree $1$ the homomorphism 
\begin{equation*}
\varphi :\mathrm{H}_{1}\left( V\right) \cong \mathbb{Z}^{k}\rightarrow 
\mathbb{Z}^{d}\cong \mathrm{H}_{1}\left( W\right) \mathrm{.}
\end{equation*}%
By the Universal Coefficient Theorem, the inclusion $V\rightarrow W$ induces
then in cohomology of degree $1$ the transpose homomorphism%
\begin{equation*}
\varphi ^{\bot }:\mathrm{H}^{1}\left( W\right) \cong \left( \mathrm{H}%
_{1}\left( W\right) \right) ^{\bot }\rightarrow \left( \mathrm{H}_{1}\left(
V\right) \right) ^{\bot }\cong \mathrm{H}^{1}\left( V\right)
\end{equation*}%
where for a free finitely-generated abelian group,%
\begin{equation*}
A^{\bot }:=\mathrm{Hom}\left( A,\mathbb{Z}\right) \text{.}
\end{equation*}

\subsection{Higher order solenoids}

Let now $A$ be a countable abelian torsion-free group. Write $A=\mathrm{co%
\mathrm{lim}}\boldsymbol{A}$ where $\boldsymbol{A}=\left( A_{n}\right) $ is
an inductive sequence of free finitely-generated groups, where $A_{n}$ has
rank $d_{n}\in \mathbb{N}$. Define $\boldsymbol{A}^{\bot }=\left(
A_{n}^{\bot }\right) $. Then $\boldsymbol{A}^{\bot }$ is a tower of free
finitely-generated groups. By the remarks above, one can then construct a
decreasing sequence $\left( W_{k}\right) $ of compact subsets of $\mathbb{R}%
^{3}$ such that:

\begin{enumerate}
\item for every $k\in \omega $, $W_{k}$ is a handlebody of genus $d_{n}$;

\item the tower $\left( \mathrm{H}_{1}\left( W_{k}\right) \right) $ is
isomorphic to $\boldsymbol{A}^{\bot }$;

\item the tower $\left( \mathrm{H}^{1}\left( S^{3}\setminus W_{k}\right)
\right) $ is isomorphic to $\boldsymbol{A}^{\bot }$;

\item the inductive sequence $\left( \mathrm{H}_{1}\left( S^{3}\setminus
W_{k}\right) \right) $ is isomorphic to $\boldsymbol{A}$.
\end{enumerate}

One then sets%
\begin{equation*}
\Sigma _{A}:=\bigcap_{k\in \omega }W_{k}\text{.}
\end{equation*}%
This is by definition a (concrete) $A$-adic solenoid, and $S^{3}\setminus
\Sigma _{A}$ is the corresponding solenoid complement. Notice that $\left(
W_{k}\right) $ provides a \emph{polyhedral resolution} for $\Sigma _{A}$ 
\cite[Theorem 9, page 65]{mardesic_shape_1982} whose bonding maps are just
the inclusion maps.

\subsection{Cohomology computation}

Fix a countable torsion-free abelian group $A$. Adopt the notation from the
previous section. Then by Steenrod (or Alexander) duality:%
\begin{equation*}
\mathrm{H}_{2}(S^{3}\setminus \Sigma _{A})\cong \mathrm{co\mathrm{lim}}_{k}{}%
\mathrm{H}_{2}(S^{3}\setminus W_{k})\cong \mathrm{co\mathrm{lim}}_{k}{}\text{%
\textrm{H}}^{2}(W_{k})=0
\end{equation*}%
and%
\begin{equation*}
\mathrm{H}_{1}(S^{3}\setminus \Sigma _{A})\cong \mathrm{co\mathrm{lim}}_{k}{}%
\mathrm{H}_{1}(S^{3}\setminus W_{k})\cong \mathrm{co\mathrm{lim}}\boldsymbol{%
A}\cong A\text{.}
\end{equation*}%
By the Universal Coefficient Theorem, this yields%
\begin{equation*}
\mathrm{H}^{2}(S^{3}\setminus \Sigma _{A};G)\cong \widetilde{\mathrm{H}}%
_{0}(\Sigma _{A};G)\cong \mathrm{Ext}(A,G)
\end{equation*}%
for every Polish abelian group $G$ and, in particular,%
\begin{equation*}
\mathrm{H}^{2}(S^{3}\setminus \Sigma _{A})\cong \mathrm{Ext}\left( A,\mathbb{%
Z}\right) \text{.}
\end{equation*}%
Similar computations show that, for $d\geq 2$,%
\begin{equation*}
\mathrm{H}_{d}(S^{3}\setminus \Sigma _{A})=0\text{ and }\mathrm{H}%
^{d+1}(S^{3}\setminus \Sigma _{A})=0\text{.}
\end{equation*}%
Thus, the \emph{Hopf Classification Theorem} produces a canonical
Borel-definable bijection%
\begin{equation*}
\lbrack S^{3}\setminus \Sigma _{A},S^{2}]\longleftrightarrow \mathrm{H}%
^{2}(S^{3}\setminus \Sigma _{A})\cong \mathrm{Ext}\left( A,\mathbb{Z}\right) 
\text{.}
\end{equation*}%
Thus, the phantom length of $[S^{3}\setminus \Sigma _{A},S^{2}]$---see
Definition \ref{Definition:phantom-length}---coincides with the phantom
length of $\mathrm{Ext}\left( A,\mathbb{Z}\right) $, which in turn is the $%
\mathbb{Z}$-projective length of $A$; see Section \ref{Section:complexity}.

\subsection{Homotopy classification}

As a consequence of the cohomology computations from the previous section we
obtain the following:

\begin{theorem}
\label{Theorem:A-adic-solenoid}Let $A$ be a countable torsion-free abelian
groups. Suppose that $\Sigma _{A}\subseteq S^{3}$ is an $A$-adic solenoids.
Then:

\begin{enumerate}
\item the classification problem for continuous maps $S^{3}\setminus \Sigma
_{A}\rightarrow S^{2}$ is equivalent to the classification problems of
extensions of $A$ by $\mathbb{Z}$;

\item every phantom map $S^{3}\setminus \Sigma _{A}\rightarrow S^{2}$ of
order at least $\alpha $ is nullhomotopic if and only if every phantom
extension of $A$ by $\mathbb{Z}$ of order at least $\alpha $ splits if and
only if $\partial _{\alpha }A$ is a free abelian group;

\item phantom maps $S^{3}\setminus \Sigma _{A}\rightarrow S^{2}$ of order $%
\alpha $ are parametrized by tails of binary sequences if and only if $%
\partial _{\alpha }A$ is the sum of a free abelian group and a finite-rank
torsion-free abelian group.
\end{enumerate}
\end{theorem}

\begin{corollary}
\label{Corollary:A-adic-solenoid}Let $\alpha <\omega _{1}$ Then:

\begin{enumerate}
\item if $\alpha $ is limit, then there exists a countable torsion-free
abelian group $A$ such that if $\Sigma _{A}\subseteq S^{3}$ is an $A$-adic
solenoid, then:

\begin{enumerate}
\item continuous maps $S^{3}\setminus \Sigma _{A}\rightarrow S^{2}$ are
classifiable up to homotopy by hereditarily countable sets of rank $1+\alpha 
$ but not of rank $1+\beta $ for any $\beta <\alpha $;

\item all maps $S^{3}\setminus \Sigma _{A}\rightarrow S^{2}$ that are
phantom of order at least $\alpha $ are trivial, but for every $\beta
<\alpha $ there exist nontrivial maps that are phantom of order at least $%
\beta $;
\end{enumerate}

\item there exists a countable torsion-free abelian group $B$ such that if $%
\Sigma _{B}\subseteq S^{3}$ is a $B$-adic solenoid, then:

\begin{enumerate}
\item continuous maps $S^{3}\setminus \Sigma _{B}\rightarrow S^{2}$ are
classifiable up to homotopy by hereditarily countable sets of $1+\alpha +1/2$
but not of rank $1+\alpha $;

\item every phantom map of order $\alpha +1$ is trivial, but the relation of
homotopy of phantom maps of order $\alpha $ is Borel bireducible with the
relation $E_{0}$ on $2^{\mathbb{N}}$ of tail equivalence of binary sequences;
\end{enumerate}

\item there exists a countable torsion-free abelian group $A$ such that if $%
\Sigma _{A}\subseteq S^{3}$ is an $A$-adic solenoid, then:

\begin{enumerate}
\item continuous maps $S^{3}\setminus \Sigma _{A}\rightarrow S^{2}$ are
classifiable up to homotopy by hereditarily countable sets of rank $1+\alpha
+1$ but not of tame rank $1+\alpha $;

\item every phantom map of order $\alpha +1$ is trivial, but the relation of
homotopy of phantom maps of order $\alpha $ is Borel bireducible with the
countable product $E_{0}^{\mathbb{N}}$ of copies of $E_{0}$.
\end{enumerate}
\end{enumerate}
\end{corollary}

\begin{proof}
This follows from Theorem \ref{Theorem:A-adic-solenoid}, Theorem \ref%
{Theorem:hereditarily-countable}, \cite[Theorem 12.16]{casarosa_phantom_2025}%
, and \cite[Theorem 6.1]{lupini_complexity_2025}.
\end{proof}

Theorem \ref{Theorem:A-adic-solenoid} and Corollary \ref%
{Corollary:A-adic-solenoid} can be seen as a higher order version of the
approach from \cite[Section 8.6]{bergfalk_definable_2024-1} to the
classification problem originally considered by Borsuk and Eilenberg in \cite%
{borsuk_uber_1936}. This problem has served as a motivation for much of the
early developments of algebraic topology and category theory as discussed in 
\cite{eilenberg_karol_1993,weibel_history_1999,kromer_tool_2007}, including 
\cite%
{eilenberg_cohomology_1940,steenrod_regular_1940,eilenberg_group_1942,eilenberg_general_1945}%
.


\bibliographystyle{amsplain}
\bibliography{cohomology}

\end{document}